\documentclass[times,sort&compress,3p]{elsarticle}
\journal{Insurance: Mathematics and Economics}
\usepackage[labelfont=bf]{caption}
\usepackage{amsmath,amsfonts,amssymb,amsthm,booktabs,color,epsfig,graphicx,url}
\usepackage{epstopdf}
\usepackage{diagbox}
\usepackage{slashbox}

\theoremstyle{plain}
\newtheorem{theorem}{Theorem}
\newtheorem{proposition}{Proposition}

\newtheorem{corollary}{Corollary}
\theoremstyle{definition}
\newtheorem{definition}{Definition}
\newtheorem{remark}{Remark}
\newtheorem{example}{Example}
\begin{document}

\begin{frontmatter}

\title{Multivariate range Value-at-Risk and covariance risk measures for elliptical and log-elliptical distributions}

\author[label1]{Baishuai  Zuo}
\author[label1]{Chuancun Yin}
\author[label2,label3]{Jing Yao\corref{mycorrespondingauthor}}
\address[label1]{School of Statistics and Data Science, Qufu Normal University, Qufu, Shandong 273165, P. R. China}
\address[label2]{Center for Financial Engineering and Department of Mathematics, Soochow University, Suzhou, P.R. China}
\address[label3]{Actuarial Mathematics and Statistics Department, Heriot-Watt University, Edinburgh, EH14 4AS, UK}

\cortext[mycorrespondingauthor]{Corresponding author. Email address: \url{j.yao@suda.edu.cn (J. Yao)}}

\begin{abstract}
In this paper, we propose the multivariate range Value-at-Risk (MRVaR) and the multivariate range covariance (MRCov) as two risk measures and explore their desirable properties in risk management. In particular, we explain that such range-based risk measures are appropriate for risk management of regulation and investment purposes. The multivariate range correlation matrix (MRCorr) is introduced accordingly. To facilitate analytical analyses, we derive explicit expressions of the MRVaR and the MRCov in the context of the multivariate (log-)elliptical distribution family. Frequently-used cases in industry, such as normal, student-$t$, logistic, Laplace, and Pearson type VII distributions, are presented with numerical examples. As an application, we propose a range-based mean-variance framework of optimal portfolio selection. We calculate the range-based efficient frontiers of the optimal portfolios based on real data of stocks' returns. Both the numerical examples and the efficient frontiers demonstrate consistences with the desirable properties of the range-based risk measures.

\end{abstract}


\end{frontmatter}

\section{Introduction}

Risk measure plays a fundamental role in many areas of insurance and
finance. Insurers and regulators rely on risk measures to calculate the risk
capital of the concerned insurance business and assess the solvency of
insurers. Investors and managers employ risk measures to evaluate the
riskiness of an investment strategy to manoeuvre optimal financial
decisions. Nowadays, risk measures such as Variance, Value-at-Risk (VaR),
tail Value-at-Risk (TVaR, also known as tail conditional expectation), and
expected shortfall (ES) are well-known to both practitioners and academics
and are ubiquitous in various quantitative analyses in actuarial risk
management and financial decision-making. Since the seminal work of Artzner
et al. (1999), abundant research works have been devoted to studying the
desirable properties of these risk measures, from both theoretical and
practical points of view. In particular, VaR is often criticised for its
non-subadditivity but enjoys elicitability, which facilitates estimation and
backtesting (Gneiting, 2011). On the contrary, TVaR is of subadditivity and
is a coherent risk measure; however, TVaR is not favoured by statisticians
as it is highly sensitive to outliers, which reduces its estimation
robustness. For a throughout discussion and review of the properties of
these risk measures, we refer to Embrechets et al. (2015), Embrechets et al.
(2018) and references therein.

To alleviate the statistical disadvantage of TVaR, Cont et al. (2010)
proposed a modified TVaR, namely the range value-at-risk (RVaR), as a bridge
between the VaR and the TVaR\bigskip . RVaR considers the average of VaR on
a given interval $[p_{1},p_{2}]$. It is defined as
\begin{equation}
\mathrm{RVaR}_{(p_{1},p_{2})}(X)=\frac{1}{p_{2}-p_{1}}\int_{p_{1}}^{p_{2}}%
\mathrm{VaR}_{t}(X)\mathrm{d}t,~0<p_{1}<p_{2}<1,  \label{(v1)}
\end{equation}%
where $X$ is the concerned risk (random variable). Note that RVaR covers
VaR, ES and TVaR as special cases, i.e.
\begin{equation*}
\lim_{p_{1}\rightarrow p_{2}}\mathrm{RVaR}_{(p_{1},p_{2})}(X)=\mathrm{VaR}%
_{p_{2}}(X),~\lim_{p_{1}\rightarrow 0}\mathrm{RVaR}_{(p_{1},p_{2})}(X)=%
\mathrm{ES}_{p_{2}}(X)~\mathrm{and}~\lim_{p_{2}\rightarrow 1}\mathrm{RVaR}%
_{(p_{1},p_{2})}(X)=\mathrm{TVaR}_{p_{1}}(X).
\end{equation*}

Many recent studies investigated financial and actuarial application
problems based on RVaR and attempted to further extend RVaR to a more
general version. For example, Embrechts et al. (2018) established an
inequality for RVaR-based risk aggregation and solved the optimal
risk-sharing problem under the RVaR. Fissler and Ziegel (2021) proposed an
elicitable triplet RVaR with two VaR components. Herrmann et al. (2020)
introduced the Geometric RVaR risk measure, which is a generalization of
RVaR and TVaR for d-dimensional distribution functions. Bairakdar et al.
(2020) defined the multivariate lower and upper orthant RVaR and study their
properties. Cai et al. (2021) defined a new multivariate conditional
Value-at-Risk (MCVaR) risk measure. As a matter of fact, RVaR can be
regarded as a conditional expectation on a truncated interval. Hence, we may
consider a multivariate conditional expectation of a truncated random vector
as a direct multivariate generalization of RVaR. Such natural generalization
is of intuition and has advantages in computation. In the literature, most
generalization is based on the tail conditional expectation (TCE). For
instance, Landsman et al. (2016a) defined a multivariate conditional tail
expectation (MTCE),
\begin{equation*}
\mathrm{MTCE}_{\boldsymbol{q}}(\mathbf{X})=\mathrm{E}\left[ \mathbf{X}|%
\mathbf{X}>\mathrm{VaR}_{\boldsymbol{q}}(\mathbf{X})\right] =\mathrm{E}[%
\mathbf{X}|X_{1}>\mathrm{VaR}_{q_{1}}(X_{1}),\cdots ,X_{n}>\mathrm{VaR}%
_{q_{n}}(X_{n})],
\end{equation*}%
where $\boldsymbol{q}=(q_{1},\cdots ,q_{n})\in (0,~1)^{n},$ $\mathbf{X}%
=(X_{1},~X_{2},\cdots ,X_{n})^{\mathrm{T}}$ is a $n\times 1$ vector of risks
with cumulative distribution function (cdf) $F_{\mathbf{X}}(\boldsymbol{x})$
and tail distribution function $\overline{F}_{\mathbf{X}}(\boldsymbol{x})$,
and $\mathrm{VaR}_{\boldsymbol{q}}(\mathbf{X})=(\mathrm{VaR}_{q_{1}}(X_{1}),~%
\mathrm{VaR}_{q_{2}}(X_{2}),\cdots ,\mathrm{VaR}_{q_{n}}(X_{n}))^{\mathrm{T}%
} $, $VaR_{q_{k}}(X_{k}),~k=1,~2,\cdots ,n,$ is the VaR of $X_{k}$.
Apparently, MTCE is reduced to tail conditional expectation (TCE) when $n=1$%
. Landsman et al. (2018) defined a multivariate tail covariance (MTCov),
\begin{equation*}
\mathrm{MTCov}_{\boldsymbol{q}}(\mathbf{X})=\mathrm{E}\left[ (\mathbf{X}-%
\mathrm{MTCE}_{\boldsymbol{q}}(\mathbf{X}))(\mathbf{X}-\mathrm{MTCE}_{%
\boldsymbol{q}}(\mathbf{X}))^{\mathrm{T}}|\mathbf{X}>\mathrm{VaR}_{%
\boldsymbol{q}}(\mathbf{X})\right] =\inf_{\boldsymbol{c}\in \mathbb{R}^{n}}%
\mathrm{E}\left[ (\mathbf{X}-\boldsymbol{c})(\mathbf{X}-\boldsymbol{c})^{%
\mathrm{T}}|\mathbf{X}>\mathrm{VaR}_{\boldsymbol{q}}(\mathbf{X})\right] .
\end{equation*}%
which is an extension of tail variance (TV) measure. Since then, several
generalizations of TCE have been proposed in the literature (see, e.g.,
Cousin and Bernardino, 2014; Cai et al., 2017; Shushi and Yao, 2020; Cai et
al., 2021; Ortega-Jim$\acute{e}$nez et al., 2021). In particular, Mousavi et
al. (2019) derived expression of MTCE for scale mixtures of skew-normal
distribution. Zuo and Yin (2021, 2022) extended those results to
location-scale mixture of elliptical and generalized skew-elliptical
distributions, respectively.

In this paper, we propose a multivariate range Value-at-Risk (MRVaR) and a
multivariate range covariance (MRCov) as risk measures using the
multivariate conditional expectation based on chosen ranges. Such
multivariate range-type (or truncated) expectation has been applied in many
areas in the literature. Arismendi and Broda (2017) derived multivariate
elliptical truncated moment generating function and pointed out that
elliptical truncated moments' expansions were useful in the design of
experiment (Thompson, 1976), robust estimation (see also Cuesta-Albertos et
al., 2008), outlier detections (Riani et al., 2009, Cerioli, 2010), robust
regression (Torti et al., 2012), robust detection (Cerioli et al. 2014)),
statistical estates' estimation (Shi et al., 2014) and risk averse selection
(Hanasusanto et al., 2015). Ogasawara (2021) derived a non-recursive formula
for various moments of the multivariate normal distribution with sectional
truncation and introduced the importance of truncated moments in biological
field. Roozegar et al. (2020) derived explicit expressions of the first two
moments for doubly truncated multivariate normal mean-variance mixture
distributions. Compared with the above mentioned studies, we make the
following contributions in this paper.

First, we formally propose the MRVaR and the MRCov as multivariate
range-based risk measures and show their desirable properties in risk
management. In particular, we point out the practical motivation to consider
such multivariate range-based risk measures. Second, as we stick to
analytical analysis, we derive explicit formulas to calculate the MRVaR and
the MRCov in the context of multivariate elliptical distribution family,
which is broadly applied in many areas. Special cases in this distribution
family and the multivariate range correlation (MRCor) are also presented.
Third, we present numerical example to further illustrate the multivariate
range-based risk measures and propose a range-based mean-variance framework
for optimal portfolio selection. We show that, by virtue of the multivariate
range-based risk measures, the range-based mean-variance framework is not
only more flexible for portfolio selection but also coherent with the
practical motivation of risk management and the desirable properties of risk
measures. As such, our MRVaR and MRCov are favourable choices from both
investors and regulator.

The rest of the paper is organized as follows. In Section 2, we give the
definitions and properties of MRVaR and MRCov risk measures, including RVaR
and RV risk measures. Section 3 introduces the elliptical and the
log-elliptical classes, and presents some special members of the elliptical
family, such as normal, student-$t$, logistic, Laplace and Pearson type VII
distributions. Section 4 derives explicit expressions of the RVaR and the
MRVaR for elliptical distributions, and also gives formulas of the RV and
the MRCov for elliptical distributions. Furthermore, those results are
extended to the class of log-elliptical distributions. Sections 5 is the
numerical illustrations. Specifically, this section compares the MRVaR and
the MRCov (MRCorr) of the several distributions, and also considers an
application of optimal portfolio selection. Finally, Section 6, is the
concluding remarks. All proofs are presented in the Appendix.

\section{Multivariate RVaR and range covariance risk measures\label{sec:2}}

Similarly with the definition of RVaR in Cont et al. (2010), we first define
the univariate RVaR and the range variance (RV) risk measures as follows.

\begin{definition}
\label{def.1} For a random variable $X$, $\mathrm{RVaR}$ risk measure is
defined by
\begin{equation*}
\mathrm{RVaR}_{(p,q)}(X)=\mathrm{E}\left[ X|\mathrm{VaR}_{p}(X)\leq X\leq
\mathrm{VaR}_{q}(X)\right] ,~p,~q\in (0,~1),~\mathrm{and}~p<q.
\end{equation*}
\end{definition}

\begin{definition}
\label{def.2} For a random variable $X$, $\mathrm{RV}$ risk measure is
defined by
\begin{align*}
\mathrm{RV}_{(p,q)}(X)=\mathrm{E}\left[(X-\mathrm{RVaR}_{(p,q)}(X))^{2}|%
\mathrm{VaR}_{p}(X)\leq X\leq \mathrm{VaR}_{q}(X)\right].
\end{align*}
\end{definition}

Note that Definition 1 is equivalent to (\ref{(v1)}), and RVaR and RV reduce
to TCE and tail variance (TV) respectively for $q\rightarrow 1$; see Furman
and Landsman (2006). Both RVaR and RV can be extended to multivariate case
straightforwardly, namely the multivariate RVaR (MRVaR) and the multivariate
range covariance (MRCov). The definitions of MRVaR and MRCov given as
follows.

\begin{definition}
\label{def.1.1} For a $n\times 1$ vector $\mathbf{X}$, $\mathrm{MRVaR}$ is
defined by
\begin{align*}
\mathrm{MRVaR}_{(\boldsymbol{p},\boldsymbol{q})}(\mathbf{X})& =\mathrm{E}%
\left[ \mathbf{X}|\mathrm{VaR}_{\boldsymbol{p}}(\mathbf{X})\leq \mathbf{X}%
\leq \mathrm{VaR}_{\boldsymbol{q}}(\mathbf{X})\right] \\
& =\mathrm{E}[\mathbf{X}|\mathrm{VaR}_{p_{1}}(X_{1})\leq X_{1}\leq \mathrm{VaR}%
_{q_{1}}(X_{1}),\cdots ,\mathrm{VaR}_{p_{n}}(X_{n})\leq X_{n}\leq \mathrm{VaR}_{q_{n}}(X_{n})],
\end{align*}%
where $\boldsymbol{p}=(p_{1},\cdots ,p_{n}),~\boldsymbol{q}=(q_{1},\cdots
,q_{n})\in (0,~1)^{n},~p_{k}<q_{k},~k=1,2,\cdots ,n.$
\end{definition}

\begin{definition}
\label{def.2.2} For a $n\times1$ vector $\mathbf{X}$, $\mathrm{MRCov}$ is
defined by
\begin{align*}
&\mathrm{MRCov}_{(\boldsymbol{p},\boldsymbol{q})}(\mathbf{X})=\mathrm{E}%
\left[(\mathbf{X}-\mathrm{MRVaR}_{(\boldsymbol{p},\boldsymbol{q})}(\mathbf{X}%
))(\mathbf{X}-\mathrm{MRVaR}_{(\boldsymbol{p},\boldsymbol{q})}(\mathbf{X}))^{%
\mathrm{T}}|\mathrm{VaR}_{\boldsymbol{p}}(\mathbf{X})\leq\mathbf{X}\leq \mathrm{VaR}_{%
\boldsymbol{q}}(\mathbf{X})\right].
\end{align*}
\end{definition}

Those range-based risk measures are advantageous in manifolds. First, they
provide straightforward and flexible extensions to many often-seen risk
measures in the literature. Taking $q_{i}=1,i=1,2,...,n$ for instance, we
immediately obtain the $\mathrm{MTCE}_{\boldsymbol{p}}$ and the $\mathrm{%
MTCov}_{\boldsymbol{p}}(\mathbf{X})$ in Landsman et al. (2016a) and Landsman
et al. (2018). Second, from a probability-actuarial point view, $\mathrm{%
MRVaR}$ and $\mathrm{MRCov}$ are also natural multivariate extensions to the
truncated expectation and covariance. In particular, comparing with other
multivariate risk measures in the literature, the calculations of proposed
measures are simpler and explicit formulas are likely reachable for given $%
\boldsymbol{p}$ and $\boldsymbol{q}$, which is important for actuarial and
financial applications. Third, the range-based risk measures are
statistically more robust in that they are less sensitive to outliers.
Moreover, these range-based risk measures are of significant practical
motivations, especially in insurance and finance; see the remarks below.

\begin{remark}
\label{re.1}In insurance and reliability engineering, range-based risk
measures arise naturally from loss models. For example, the deductible and
the policy limit in insurance contracts are of typical range-based forms
(Klugman et al. 1998). From a data analysis point of view, range-based risk
measures also provide convenient mathematical expressions for survival data
and failure time data, such as time to onset of a disease, length of stay in
an automobile insurance plan (company), money paid for hospitalization by
health insurance, etc. Such data is very common in insurance and reliability
and can be only interval-observed in many cases.

In finance and risk management, regulators and decision-makers usually
concern the \textquotedblleft tails\textquotedblright\ of risks because
firms' financial situations are severely impacted by extreme events, such as
businesses exposed to property-catastrophe risks for instance. In fact, this
is exactly the motivation that boosts many risk models in the first place.
On the other hand, from the investors' point of view, it is more appropriate
to study the non-tail part of the risk, which stands for the firm's ordinary
business and the shareholders' interests. Clearly, our range-based risk
measures are applicable for both.
\end{remark}

The above remarks also suggests feasible choices of the ranges. Regulatory
constrains, investment target and risk management requirement (e.g. credit
rating level) can be considered in determining desirable $\boldsymbol{p}$
and $\boldsymbol{q}$. To formulate the proposed range-based risk measures,
we further study the following properties of the MRVaR and the MRCov.

\begin{proposition}
\label{pro.1} For any $n\times1$ random vectors $\mathbf{X}%
=(X_{1},X_{2},\cdots,X_{n})^{\mathrm{T}}$ and $\mathbf{Y}=(Y_{1},Y_{2},%
\cdots,Y_{n})^{\mathrm{T}}$, the $\mathrm{MRVaR}$ risk measure has the
following properties:\newline
(i) (Positive homogeneity) For any positive constant $c$, we have
\begin{equation*}
\mathrm{MRVaR}_{(\boldsymbol{p},\boldsymbol{q})}(c\mathbf{X})=c\mathrm{MRVaR}%
_{(\boldsymbol{p},\boldsymbol{q})}(\mathbf{X});
\end{equation*}
(ii) (Translation invariance) For any vector of constants $\boldsymbol{\gamma%
}\in \mathbb{R}^{n}$
\begin{equation*}
\mathrm{MRVaR}_{(\boldsymbol{p},\boldsymbol{q})}(\mathbf{X}+\boldsymbol{%
\gamma})=\mathrm{MRVaR}_{(\boldsymbol{p},\boldsymbol{q})}(\mathbf{X})+%
\boldsymbol{\gamma};
\end{equation*}
(iii) (Independency of risks) When $\mathbf{X}$ has independent components,
we have
\begin{equation*}
\mathrm{MRVaR}_{(\boldsymbol{p},\boldsymbol{q})}(\mathbf{X})=(\mathrm{RVaR}%
_{(p_{1},q_{1})}(X_{1}),\mathrm{RVaR}_{(p_{2},q_{2})}(X_{2}),\cdots,\mathrm{%
RVaR}_{(p_{n},q_{n})}(X_{n}))^\mathrm{T};
\end{equation*}
(iv) (Monotonicity) If $\mathbf{Y}\overset{a.s.}{\geq}\mathbf{X}$, then
\begin{equation*}
\mathrm{MRVaR}_{(\boldsymbol{p},\boldsymbol{q})}(\mathbf{Y}-\mathbf{X})\geq
\boldsymbol{0},
\end{equation*}
where $\boldsymbol{0}$ is vector of $n$ zeros.
\end{proposition}

\begin{proposition}
\label{pro.2} For any $n\times1$ random vector $\mathbf{X}%
=(X_{1},X_{2},\cdots,X_{n})^{\mathrm{T}}$, the $\mathrm{MRCov}$ risk measure
has following properties:\newline
(i)(Standardization) For any vector of constants $\boldsymbol{\gamma}\in
\mathbb{R}^{n}$
\begin{equation*}
\mathrm{MRCov}_{(\boldsymbol{p},\boldsymbol{q})}(\boldsymbol{\gamma})=%
\boldsymbol{0}_{n\times n},
\end{equation*}
where $\boldsymbol{0}_{n\times n}$ is an $n\times n$ matrix, whose
components are zeros; \newline
(ii) (Positive homogeneity) For any positive constant $c$, we have
\begin{equation*}
\mathrm{MRCov}_{(\boldsymbol{p},\boldsymbol{q})}(c\mathbf{X})=c^{2}\mathrm{%
MRCov}_{(\boldsymbol{p},\boldsymbol{q})}(\mathbf{X});
\end{equation*}
(iii) (Translation invariance) For any vector of constants $\boldsymbol{%
\gamma}\in \mathbb{R}^{n}$
\begin{equation*}
\mathrm{MRCov}_{(\boldsymbol{p},\boldsymbol{q})}(\mathbf{X}+\boldsymbol{%
\gamma})=\mathrm{MRCov}_{(\boldsymbol{p},\boldsymbol{q})}(\mathbf{X});
\end{equation*}
(iv) (Independency of risks) When $\mathbf{X}$ has independent components,
we have
\begin{equation*}
\mathrm{MRCov}_{(\boldsymbol{p},\boldsymbol{q})}(\mathbf{X})=\mathrm{diag}(%
\mathrm{RV}_{(p_{1},q_{1})}(X_{1}),\mathrm{RV}_{(p_{2},q_{2})}(X_{2}),\cdots,%
\mathrm{RV}_{(p_{n},q_{n})}(X_{n})),
\end{equation*}
where $\mathrm{diag}()$ is diagonal matrix.
\end{proposition}

From matrix MRCov, we can definite the multivariate range correlation
(MRCorr) matrix:
\begin{equation}  \label{(v2)}
\mathrm{MRCorr}_{(\boldsymbol{p},\boldsymbol{q})}(\mathbf{X})=\bigg(\frac{%
\mathrm{MRCov}_{(\boldsymbol{p},\boldsymbol{q})}(\mathbf{X})_{ij}}{\sqrt{%
\mathrm{MRCov}_{(\boldsymbol{p},\boldsymbol{q})}(\mathbf{X})_{ii}}\sqrt{%
\mathrm{MRCov}_{(\boldsymbol{p},\boldsymbol{q})}(\mathbf{X})_{jj}}}\bigg)%
_{ij=1,\cdots ,n}.
\end{equation}

It is our primary intention to derive explicit expressions for relevant
analysis in this paper. To this end, we work in the context of multivariate
elliptical and log-elliptical distributions. We introduce the two
probability distribution families in next section.

\section{Multivariate Elliptical and Log-elliptical Distributions}

\subsection{Multivariate elliptical distributions}

We say a $n$-dimensional random vector $\mathbf{X}=(X_{1},\cdots ,X_{n})^{%
\mathrm{T}}$ follows a multivariate elliptical distribution if its
characteristic function has the form
\begin{equation*}
\mathrm{E}[\exp (i\mathbf{t}^{T}\mathbf{X})]=e^{i\mathbf{t}^{T}{\boldsymbol{%
\mu }}}\psi_{n} \left( \frac{1}{2}\mathbf{t}^{T}\mathbf{\Sigma }\mathbf{t}%
\right)
\end{equation*}%
for all $\mathbf{t}\in \mathbb{R}^{n}$, denoted $\mathbf{X}\sim E_{n}({%
\boldsymbol{\mu }},\mathbf{\Sigma },\psi_{n} )$, where $\psi_{n} $ is called
the characteristic generator with $\psi_{n} (0)=1$ (Fang et al., 1990). If
probability density function (pdf) of $\mathbf{X}$ exists, it is of the form
below,
\begin{equation*}
f_{\boldsymbol{X}}(\boldsymbol{x}):=\frac{c_{n}}{\sqrt{|\boldsymbol{\Sigma }|%
}}g_{n}\left\{ \frac{1}{2}(\boldsymbol{x}-\boldsymbol{\mu })^{T}\mathbf{%
\Sigma }^{-1}(\boldsymbol{x}-\boldsymbol{\mu })\right\} ,~\boldsymbol{x}\in
\mathbb{R}^{n},
\end{equation*}%
where $\boldsymbol{\mu }$ is a $n\times 1$ location vector, $\mathbf{\Sigma }
$ is a $n\times n$ scale matrix, and $g_{n}(u)$, $u\geq 0$, is called the
density generator of $\mathbf{X}$. We then denote it as $\mathbf{X}\sim
E_{n}(\boldsymbol{\mu },\boldsymbol{\Sigma },g_{n})$. The density generator $%
g_{n}$ must satisfy the condition
\begin{equation}
\int_{0}^{\infty }s^{n/2-1}g_{n}(s)\mathrm{d}s<\infty ,  \label{(v3)}
\end{equation}%
such that $f_{\boldsymbol{X}}(\boldsymbol{x})$ is a density function, where
the normalizing constant $c_{n}$ is
\begin{equation*}
c_{n}=\frac{\Gamma (n/2)}{(2\pi )^{n/2}}\left[ \int_{0}^{\infty
}s^{n/2-1}g_{n}(s)\mathrm{d}s\right] ^{-1}.
\end{equation*}%
Suppose $\mathbf{A}$ is a $k\times n$ matrix ($k\leq n$), and $\mathbf{b}$
is a $k\times 1$ vector. Then
\begin{equation}
\mathbf{AX+b}\sim E_{k}\left( \boldsymbol{A\mu +b},~\boldsymbol{A\Sigma A^{T}%
},~g_{k,n}\right) ,  \label{(v4)}
\end{equation}%
where $g_{k,n}(u)=\int_{0}^{\infty }s^{\frac{n-k}{2}-1}g_{n}(s+u)\mathrm{d}s$%
.

To facilitate the calculation of the proposed risk measures, we further
define two associated random vectors of $\mathbf{X}$. For a density
generator $g_{n}(s)$, we define its cumulative generator $\overline{G}%
_{n}(u) $ and $\overline{\mathcal{G}}_{n}(u)$ as
\begin{equation*}
\overline{G}_{n}(u)=\int_{u}^{\infty }{g}_{n}(v)\mathrm{d}v~~\mathrm{and}~~%
\overline{\mathcal{G}}_{n}(u)=\int_{u}^{\infty }{\overline{G}}_{n}(v)\mathrm{%
d}v,
\end{equation*}%
and the normalizing constants are, respectively, written as
\begin{equation*}
c_{n}^{\ast }=\frac{\Gamma (n/2)}{(2\pi )^{n/2}}\left[ \int_{0}^{\infty
}s^{n/2-1}\overline{G}_{n}(s)\mathrm{d}s\right] ^{-1}~~\mathrm{and}%
~~c_{n}^{\ast \ast }=\frac{\Gamma (n/2)}{(2\pi )^{n/2}}\left[
\int_{0}^{\infty }s^{n/2-1}\overline{\mathcal{G}}_{n}(s)\mathrm{d}s\right]
^{-1}.
\end{equation*}%
It can be verified that both $\overline{G}_{n}(u)$ and $\overline{\mathcal{G}%
}_{n}(u)$ satisfy the required conditions
\begin{equation}  \label{(v5)}
\int_{0}^{\infty }s^{n/2-1}\overline{G}_{n}(s)\mathrm{d}s<+\infty
\end{equation}%
and
\begin{equation}  \label{(v6)}
\int_{0}^{\infty }s^{n/2-1}\overline{\mathcal{G}}_{n}(s)\mathrm{d}s<+\infty .
\end{equation}%
Then $\mathbf{X}^{\ast }\sim E_{n}(\boldsymbol{\mu },~\boldsymbol{\Sigma },~%
\overline{G}_{n})$ and $\mathbf{X}^{\ast \ast }\sim E_{n}(\boldsymbol{\mu },~%
\boldsymbol{\Sigma },~\overline{\mathcal{G}}_{n})$ (Zuo et al., 2021) are
elliptical random vectors with density generators $\overline{G}_{n}(u)$ and $%
\overline{\mathcal{G}}_{n}(u)$, and their density functions are
\begin{equation*}
f_{\boldsymbol{X}^{\ast }}(\boldsymbol{x})=\frac{c_{n}^{\ast }}{\sqrt{|%
\boldsymbol{\Sigma }|}}\overline{G}_{n}\left\{ \frac{1}{2}(\boldsymbol{x}-%
\boldsymbol{\mu })^{T}\mathbf{\Sigma }^{-1}(\boldsymbol{x}-\boldsymbol{\mu }%
)\right\} ,~\boldsymbol{x}\in \mathbb{R}^{n},
\end{equation*}%
\begin{equation*}
f_{\boldsymbol{X}^{\ast \ast }}(\boldsymbol{x})=\frac{c_{n}^{\ast \ast }}{%
\sqrt{|\boldsymbol{\Sigma }|}}\overline{\mathcal{G}}_{n}\left\{ \frac{1}{2}(%
\boldsymbol{x}-\boldsymbol{\mu })^{T}\mathbf{\Sigma }^{-1}(\boldsymbol{x}-%
\boldsymbol{\mu })\right\} ,~\boldsymbol{x}\in \mathbb{R}^{n}.
\end{equation*}

The following examples are special members of the elliptical family that
have been frequently-used in modelling risks in various areas.

\begin{example}
\label{exa.1} (Multivariate normal distribution). Suppose that $\mathbf{X}%
\sim N_{n}\left(\boldsymbol{\mu},~\boldsymbol{\Sigma}\right)$. In this case,
the density generators $g_{n}(u)$, $\overline{G}_{n}(u)$ and $\overline{%
\mathcal{G}}_{n}(u)$ are expressed:
\begin{align*}
g(u)=\overline{G}(u)=\overline{\mathcal{G}}(u)=\exp\{-u\},
\end{align*}
and the normalizing constants are written:
\begin{align*}
c_{n}=c_{n}^{\ast}=c_{n}^{\ast\ast}=(2\pi)^{-\frac{n}{2}}.
\end{align*}
\end{example}

\begin{example}
\label{exa.2} (Multivariate student-$t$ distribution). Suppose that
\begin{align}  \label{(v7)}
\mathbf{X}\sim St_{n}\left(\boldsymbol{\mu},~\boldsymbol{\Sigma},~m\right).
\end{align}
In this case, the density generators $g_{n}(u)$, $\overline{G}_{n}(u)$ and $%
\overline{\mathcal{G}}_{n}(u)$ are expressed (for details, see Zuo et al.,
2021):
\begin{align*}
g_{n}(u)=\left(1+\frac{2u}{m}\right)^{-(m+n)/2}, ~~ \overline{G}_{n}(u)=%
\frac{m}{m+n-2}\left(1+\frac{2u}{m}\right)^{-(m+n-2)/2}
\end{align*}
and
\begin{align*}
\overline{\mathcal{G}}_{n}(u)=\frac{m}{m+n-2}\frac{m}{m+n-4}\left(1+\frac{2u%
}{m}\right)^{-(m+n-4)/2}.
\end{align*}
The normalizing constants are written:
\begin{align*}
c_{n}=\frac{\Gamma\left((m+n)/2\right)}{\Gamma(m/2)(m\pi)^{\frac{n}{2}}},
\end{align*}
\begin{align*}
c_{n}^{\ast}&=\frac{(m+n-2)\Gamma(n/2)}{(m\pi)^{n/2}mB(\frac{n}{2},~\frac{m-2%
}{2})},~\mathrm{if}~m>2
\end{align*}
and
\begin{align*}
c_{n}^{\ast\ast}&=\frac{(m+n-2)(m+n-4)\Gamma(n/2)}{(m\pi)^{n/2}m^{2}B(\frac{n%
}{2},~\frac{m-4}{2})},~\mathrm{if}~m>4,
\end{align*}
where $\Gamma(\cdot)$ and $B(\cdot,\cdot)$ are Gamma function and Beta
function, respectively.
\end{example}

\begin{example}
\label{exa.3} (Multivariate logistic distribution). Suppose that $\mathbf{X}%
\sim Lo_{n}\left( \boldsymbol{\mu },~\boldsymbol{\Sigma }\right) $. In this
case, the density generators $g_{n}(u)$, $\overline{G}_{n}(u)$ and $%
\overline{\mathcal{G}}_{n}(u)$ are expressed (see also Zuo et al., 2021):
\begin{equation*}
g_{n}(u)=\frac{\exp (-u)}{[1+\exp (-u)]^{2}},~~\overline{G}_{n}(u)=\frac{%
\exp (-u)}{1+\exp (-u)}~~\mathrm{and}~~\overline{\mathcal{G}}_{n}(u)=\ln %
\left[ 1+\exp (-u)\right] .
\end{equation*}%
The normalizing constants are written:
\begin{equation*}
c_{n}=\frac{1}{(2\pi )^{n/2}\Psi _{2}^{\ast }(-1,\frac{n}{2},1)}%
,~~c_{n}^{\ast }=\frac{1}{(2\pi )^{n/2}\Psi _{1}^{\ast }(-1,\frac{n}{2},1)}~~%
\mathrm{and}~~c_{n}^{\ast \ast }=\frac{1}{(2\pi )^{n/2}\Psi _{1}^{\ast }(-1,%
\frac{n}{2}+1,1)}.
\end{equation*}
\end{example}

\begin{remark}
\label{re.2} Here $\Psi_{\kappa}^{\ast}(z,s,a)$ is the generalized
Hurwitz-Lerch zeta function defined by (see Lin et al., 2006)
\begin{equation*}
\Psi_{\kappa}^{\ast}(z,s,a)=\frac{1}{\Gamma(\kappa)}\sum_{n=0}^{\infty}\frac{%
\Gamma(\kappa+n)}{n!}\frac{z^{n}}{(n+a)^{s}},
\end{equation*}
which has an integral representation
\begin{equation*}
\Psi_{\kappa}^{\ast}(z,s,a)=\frac{1}{\Gamma(s)}\int_{0}^{\infty}\frac{%
t^{s-1}e^{-at}}{(1-ze^{-t})^{\kappa}}\mathrm{d}t,
\end{equation*}
where $\mathcal{R}(a)>0$, $\mathcal{R}(s)>0$ when $|z|\leq1~(z\neq1)$, $%
\mathcal{R}(s)>1$ when $z=1$.
\end{remark}

\begin{example}
\label{exa.4} (Multivariate Laplace distribution). Suppose that $\mathbf{X}%
\sim La_{n}\left( \boldsymbol{\mu },~\boldsymbol{\Sigma }\right) $. In this
case, the density generators $g_{n}(u)$, $\overline{G}_{n}(u)$ and $%
\overline{\mathcal{G}}_{n}(u)$ are expressed (see also Zuo et al., 2021):
\begin{equation*}
g_{n}(u)=\exp (-\sqrt{2u}),~~\overline{G}_{n}(u)=(1+\sqrt{2u})\exp (-\sqrt{2u%
})
\end{equation*}%
and
\begin{equation*}
\overline{\mathcal{G}}_{n}(u)=(3+2u+3\sqrt{2u})\exp (-\sqrt{2u}).
\end{equation*}%
The normalizing constants are written:
\begin{equation*}
c_{n}=\frac{\Gamma (n/2)}{2\pi ^{n/2}\Gamma (n)},~~c_{n}^{\ast }=\frac{%
n\Gamma (n/2)}{2\pi ^{n/2}\Gamma (n+2)}~~\mathrm{and}~~c_{n}^{\ast \ast }=%
\frac{n(n+2)\Gamma (n/2)}{2\pi ^{n/2}\Gamma (n+4)}.
\end{equation*}
\end{example}

\begin{example}
\label{exa.5} (Multivariate Pearson type VII distribution). Suppose that $%
\mathbf{X}\sim PVII_{n}\left(\boldsymbol{\mu},~\boldsymbol{\Sigma}%
,~t\right). $ In this case, the density generators $g_{n}(u)$, $\overline{G}%
_{n}(u)$ and $\overline{\mathcal{G}}_{n}(u)$ are expressed:
\begin{align*}
g_{n}(u)=(1+2u)^{-t}, ~~ \overline{G}_{n}(u)=\frac{1}{2(t-1)}(1+2u)^{-(t-1)}
\end{align*}
and
\begin{align*}
\overline{\mathcal{G}}_{n}(u)=\frac{1}{4(t-1)(t-2)}(1+2u)^{-(t-2)}.
\end{align*}
The normalizing constants are written:
\begin{align*}
c_{n}=\frac{\Gamma\left(t\right)}{\Gamma(t-n/2)\pi^{\frac{n}{2}}},~t>\frac{n%
}{2},
\end{align*}
\begin{align*}
c_{n}^{\ast}=\frac{\Gamma(\frac{n}{2})2(t-1)}{\pi^{\frac{n}{2}}B(\frac{n}{2}%
,t-1-\frac{n}{2})},~t>1+\frac{n}{2}
\end{align*}
and
\begin{align*}
c_{n}^{\ast\ast}=\frac{\Gamma(\frac{n}{2})4(t-1)(t-2)}{\pi^{\frac{n}{2}}B(%
\frac{n}{2},t-2-\frac{n}{2})},~t>2+\frac{n}{2}.
\end{align*}
\end{example}

\begin{remark}
\label{re.3} Multivariate Pearson type VII distribution is related to the
multivariate Student-$t$ by the transformation $\mathbf{X}=[\mathbf{Y}-(1-%
\sqrt{m}\boldsymbol{\mu})]/\sqrt{m}$, where $\mathbf{Y}$ is as in (\ref{(v7)}%
) with $m=2t-n$ (see Zografos, 2008).
\end{remark}

\subsection{Multivariate log-elliptical distributions}

For any $n$-dimensional vector $\mathbf{Z}=(Z_{1},Z_{2},\cdots ,Z_{n})^{T}$
with positive components $Z_{k}$, $k=1,2,\cdots ,n,$ we write $\mathrm{ln}%
\mathbf{Z}=(\mathrm{ln}Z_{1},\mathrm{ln}Z_{2},\cdots ,\mathrm{\ ln}%
Z_{n})^{T} $. If random vector $\mathrm{ln}\mathbf{Z}\sim E_{n}(\boldsymbol{%
\mu },\mathbf{\ \Sigma },\psi_{n} )$, we will say that the random vector $%
\mathbf{Z}$ has an $n$- variate log-elliptical distribution, denoted by $%
\mathbf{Z}\sim LE_{n}(\boldsymbol{\mu },\mathbf{\Sigma },\psi_{n} )$.
Further, if its probability density function (pdf) exists, the form will be
(see Valdez et al., 2009)
\begin{equation*}
f_{\boldsymbol{Z}}(\boldsymbol{z}):=\frac{c_{n}}{\sqrt{|\boldsymbol{\Sigma }|%
}}\left( \prod_{k=1}^{n}z_{k}^{-1}\right) g_{n}\left\{ \frac{1}{2}(\mathrm{ln%
}\boldsymbol{z}-\boldsymbol{\mu })^{T}\mathbf{\Sigma }^{-1}(\mathrm{ln}%
\boldsymbol{z}-\boldsymbol{\mu })\right\} ,~z_{k}>0,~\forall k\in
\{1,2,\cdots ,n\},
\end{equation*}%
where $\mathrm{ln}\boldsymbol{z}=(\mathrm{ln}z_{1},\mathrm{ln}z_{2},\cdots ,%
\mathrm{ln}z_{n})^{\mathrm{T}}$. In this case we denote that $\mathbf{Z}\sim
LE_{n}(\boldsymbol{\mu },\mathbf{\Sigma },g_{n})$.

\section{Main Results}

\subsection{RVaR and MRVaR of elliptical distribution}

For an elliptical random vector $\mathbf{X}\sim E_{n}\left( \boldsymbol{\mu }%
,~\boldsymbol{\Sigma },~g_{n}\right) $ with density function $f_{\boldsymbol{%
X}}(\boldsymbol{x})$. It is evident that $\mathbf{Y}=\mathbf{\Sigma }^{-%
\frac{1}{2}}(\mathbf{X}-\boldsymbol{\mu })\sim E_{n}\left( \boldsymbol{0},~%
\boldsymbol{I_{n}},~g_{n}\right) .$ We employ the following notations in the
sequel of the paper.
\begin{equation*}
\boldsymbol{\eta _{v}}=\left( \eta _{\boldsymbol{v},1},~\eta _{\boldsymbol{v}%
,2},\cdots ,\eta _{\boldsymbol{v},n}\right) ^{\mathrm{\mathrm{T}}}=\mathbf{%
\Sigma }^{-\frac{1}{2}}(\mathrm{VaR}_{\boldsymbol{v}}-\boldsymbol{\mu }),
\end{equation*}%
where
\begin{equation*}
\boldsymbol{\eta }_{\boldsymbol{v},-k}=\left( \eta _{\boldsymbol{v},1},~\eta
_{\boldsymbol{v},2},\cdots ,\eta _{\boldsymbol{v},k-1},~\eta _{\boldsymbol{v}%
,k+1},\cdots ,\eta _{\boldsymbol{v},n}\right) ^{\mathrm{T}}
\end{equation*}%
and
\begin{equation*}
\boldsymbol{\eta }_{\boldsymbol{v},-k,j}=\left( \eta _{\boldsymbol{v}%
,1},\cdots ,\eta _{\boldsymbol{v},k-1},~\eta _{\boldsymbol{v},k+1},\cdots
,\eta _{\boldsymbol{v},j-1},~\eta _{\boldsymbol{v},j+1},\cdots ,\eta _{%
\boldsymbol{v},n}\right) ^{\mathrm{T}},~\boldsymbol{v}\in \{\boldsymbol{p},%
\boldsymbol{q}\}.
\end{equation*}%
A truncated distribution function is defined as
\begin{equation*}
F_{\mathbf{Z}}(\boldsymbol{a},\boldsymbol{b})=\int_{\boldsymbol{a}}^{%
\boldsymbol{b}}f_{\mathbf{Z}}(\boldsymbol{z})\mathrm{d}\boldsymbol{z},
\end{equation*}%
where $f_{\mathbf{Z}}(\boldsymbol{z})$ is pdf of random vector $\mathbf{Z}$.
Note that that $F_{\mathbf{Z}}(\boldsymbol{a},\boldsymbol{b})=F_{\mathbf{Z}}(%
\boldsymbol{b})$ as $\boldsymbol{a}\rightarrow \mathbf{-\infty }$, and $F_{%
\mathbf{Z}}(\boldsymbol{a},\boldsymbol{b})=\overline{F}_{\mathbf{Z}}(%
\boldsymbol{b})$ as $\boldsymbol{b}\rightarrow \mathbf{+\infty }$.\newline
For multivariate elliptical distributions, we also define shifted cumulative
generator
\begin{equation*}
\overline{G}_{n-1,\boldsymbol{v},k}^{\ast }(u)=\overline{G}_{n}\left( u+%
\frac{1}{2}\eta _{\boldsymbol{v},k}^{2}\right) ,~k=1,2,\cdots ,n,~%
\boldsymbol{v}\in \{\boldsymbol{p},\boldsymbol{q}\},
\end{equation*}%
and normalizing constant
\begin{equation*}
c_{n-1,\boldsymbol{v},k}^{\ast }=\frac{\Gamma (\frac{n-1}{2})}{(2\pi
)^{(n-1)/2}}\left[ \int_{0}^{\infty }s^{(n-3)/2}\overline{G}_{n-1,%
\boldsymbol{v},k}^{\ast }(s)\mathrm{d}s\right] ^{-1}.
\end{equation*}

\begin{theorem}
\label{th.1.1} Let $X\sim E_{1}(\mu ,~\sigma ^{2},~g_{1})$ be an elliptical
random variable with density generator $g_{1}$, scale parameter $\sigma $
and finite expectation $\mu $. We have
\begin{equation*}
\mathrm{RVaR}_{(p,q)}(X)=\mu +\sigma \frac{\delta _{p,q}}{F_{Y}(\eta
_{p},\eta _{q})},
\end{equation*}%
where
\begin{equation}  \label{(v8)}
\delta _{p,q}=c_{1}\left[ \overline{G}_{1}\left( \frac{1}{2}\eta
_{p}^{2}\right) -\overline{G}_{1}\left( \frac{1}{2}\eta _{q}^{2}\right) %
\right] ,
\end{equation}%
and $Y\sim E_{1}\left( 0,~1,~g_{1}\right) $.
\end{theorem}

\begin{theorem}
\label{th.1} Let $\mathbf{X}\sim E_{n}(\boldsymbol{\mu },~\mathbf{\Sigma }%
,~g_{n})$ $(n\geq 2)$ be an $n$-dimensional elliptical random vector with
density generator $g_{n}$, positive defined scale matrix $\mathbf{\Sigma }$
and finite expectation $\boldsymbol{\mu }$. Then,
\begin{equation*}
\mathrm{MRVaR}_{(\boldsymbol{p,q})}(\mathbf{X})=\boldsymbol{\mu }+\mathbf{%
\Sigma }^{\frac{1}{2}}\frac{\boldsymbol{\delta _{\boldsymbol{p,q}}}}{F_{%
\mathbf{Y}}(\boldsymbol{\eta }_{\boldsymbol{p}},\boldsymbol{\eta }_{%
\boldsymbol{q}})},
\end{equation*}%
where
\begin{equation*}
\boldsymbol{\delta _{\boldsymbol{p,q}}}=\left( \delta _{1,\boldsymbol{p,q}%
},~\delta _{2,\boldsymbol{p,q}},\cdots ,\delta _{n,\boldsymbol{p,q}}\right)
^{T}
\end{equation*}%
with
\begin{equation*}
\delta _{k,\boldsymbol{p,q}}=\frac{c_{n}}{c_{n-1,\boldsymbol{p},k}^{\ast }}%
F_{\mathbf{Y}_{\boldsymbol{p},-k}}(\boldsymbol{\eta }_{\boldsymbol{p},-k},%
\boldsymbol{\eta }_{\boldsymbol{q},-k})-\frac{c_{n}}{c_{n-1,\boldsymbol{q}%
,k}^{\ast }}F_{\mathbf{Y}_{\boldsymbol{q},-k}}(\boldsymbol{\eta }_{%
\boldsymbol{p},-k},\boldsymbol{\eta }_{\boldsymbol{q},-k}),~k=1,~2,\cdots ,n,
\end{equation*}%
$\mathbf{Y}\sim E_{n}\left( \boldsymbol{0},~\boldsymbol{I_{n}},~g_{n}\right)
$, $\mathbf{Y}_{\boldsymbol{v},-k}\sim E_{n-1}\left( \boldsymbol{0},~%
\boldsymbol{I_{n-1}},~\overline{G}_{n-1,\boldsymbol{v},k}^{\ast }\right) $.
\end{theorem}

\begin{remark}
\label{re.4} Note that $\boldsymbol{q}\rightarrow (1,\cdots ,1)^{T}$ in
Theorem \ref{th.1}, one gets the formula of Theorem 1 in Landsman et al.
(2018), which is generalization of Theorem 1 in Landsman et al. (2016a); in
the one-dimensional case, we obtain result of Theorem 1 in Landsman and
Valdez (2003), which is the generalization of Theorem 1 in Jiang et al.
(2016).
\end{remark}

By replacing the density generators $g_{n}(u)$, $\overline{G}_{n}(u)$ and $%
\overline{\mathcal{G}}_{n}(u)$ of Examples \ref{exa.1}-\ref{exa.5} into
Theorem \ref{th.1}, and doing some algebraic operations, we can obtain
following corollaries.

\begin{corollary}
\label{co.1} Let $\mathbf{X}\sim N_{n}(\boldsymbol{\mu},~\mathbf{\Sigma})$ $%
(n\geq2)$. Then
\begin{align*}
&\mathrm{MRVaR}_{(\boldsymbol{p,q})}(\mathbf{X})=\boldsymbol{\mu}+\mathbf{%
\Sigma}^{\frac{1}{2}}\frac{\boldsymbol{\delta_{\boldsymbol{p,q}}}}{F_{%
\mathbf{Y}}(\boldsymbol{\eta}_{\boldsymbol{p}},\boldsymbol{\eta}_{%
\boldsymbol{q}})},
\end{align*}
where
\begin{equation*}
\boldsymbol{\delta_{\boldsymbol{p,q}}}=\left(\delta_{1,\boldsymbol{p,q}%
},~\delta_{2,\boldsymbol{p,q}},\cdots,\delta_{n,\boldsymbol{p,q}}\right)^{T}
\end{equation*}
with
\begin{align*}
\delta_{k,\boldsymbol{p,q}}&=\left[\phi(\eta_{\boldsymbol{p},k})-\phi(\eta_{%
\boldsymbol{q},k})\right]F_{\mathbf{Y}_{-k}}(\boldsymbol{\eta}_{\boldsymbol{p%
},-k},\boldsymbol{\eta}_{\boldsymbol{q},-k}),~k=1,~2,\cdots,n,
\end{align*}
$\mathbf{Y}\sim N_{n}\left(\boldsymbol{0},~\boldsymbol{I_{n}}\right)$, $%
\mathbf{Y}_{-k}\sim N_{n-1}\left(\boldsymbol{0},~\boldsymbol{I_{n-1}}\right)$%
, and $\phi(\cdot)$ is pdf of $1$-dimensional standard normal distribution.
\end{corollary}

\begin{corollary}
\label{co.2} Let $\mathbf{X}\sim St_{n}(\boldsymbol{\mu},~\mathbf{\Sigma}%
,~m) $ $(n\geq2)$. Then
\begin{align*}
&\mathrm{MRVaR}_{(\boldsymbol{p,q})}(\mathbf{X})=\boldsymbol{\mu}+\mathbf{%
\Sigma}^{\frac{1}{2}}\frac{\boldsymbol{\delta_{\boldsymbol{p,q}}}}{F_{%
\mathbf{Y}}(\boldsymbol{\eta}_{\boldsymbol{p}},\boldsymbol{\eta}_{%
\boldsymbol{q}})},
\end{align*}
where
\begin{equation*}
\boldsymbol{\delta_{\boldsymbol{p,q}}}=\left(\delta_{1,\boldsymbol{p,q}%
},~\delta_{2,\boldsymbol{p,q}},\cdots,\delta_{n,\boldsymbol{p,q}}\right)^{T}
\end{equation*}
with
\begin{align*}
&\delta_{k,\boldsymbol{p,q}}=\frac{c_{n}}{c_{n-1,\boldsymbol{p},k}^{\ast}}F_{%
\mathbf{Y}_{\boldsymbol{p},-k}}(\boldsymbol{\eta}_{\boldsymbol{p},-k},%
\boldsymbol{\eta}_{\boldsymbol{q},-k})-\frac{c_{n}}{c_{n-1,\boldsymbol{q}%
,k}^{\ast}}F_{\mathbf{Y}_{\boldsymbol{q},-k}}(\boldsymbol{\eta}_{\boldsymbol{%
p},-k},\boldsymbol{\eta}_{\boldsymbol{q},-k}),~k=1,~2,\cdots,n,
\end{align*}
$\mathbf{Y}\sim St_{n}\left(\boldsymbol{0},~\boldsymbol{I_{n}},~m\right)$, $%
\mathbf{Y}_{\boldsymbol{v},-k}\sim St_{n-1}\left(\boldsymbol{0},~\Delta_{%
\boldsymbol{v},k},~m-1\right),$
\begin{equation*}
\Delta_{\boldsymbol{v},k}=\left[\frac{m\left(1+\eta_{\boldsymbol{v}%
,k}^{2}/m\right)}{m-1}\right]\mathbf{I}_{n-1}
\end{equation*}
and
\begin{equation*}
c_{n-1,\boldsymbol{v},k}^{\ast}=\frac{\Gamma\left(\frac{m+n-2}{2}%
\right)(m+n-2)}{\Gamma\left(\frac{m-1}{2}\right)[(m-1)\pi]^{(n-1)/2}m}%
\left(1+\frac{\eta_{\boldsymbol{v},k}^{2}}{m}\right)^{(m+n-2)/2},~%
\boldsymbol{v}\in\{\boldsymbol{p},\boldsymbol{q}\}.
\end{equation*}
\end{corollary}

We can further simplify
\begin{equation*}
\frac{c_{n}}{c_{\boldsymbol{v},k}^{\ast}}=\frac{\Gamma\left(\frac{m-1}{2}%
\right)\left(\frac{m-1}{m}\right)^{(n-1)/2}}{\Gamma\left(\frac{m}{2}\right)%
\sqrt{\pi/m}}\left(1+\frac{\eta_{\boldsymbol{v},k}^{2}}{m}%
\right)^{-(m+n-2)/2},~\boldsymbol{v}\in\{\boldsymbol{p},\boldsymbol{q}\}.
\end{equation*}

\begin{corollary}
\label{co.3} Let $\mathbf{X}\sim Lo_{n}(\boldsymbol{\mu},~\mathbf{\Sigma})$ $%
(n\geq2)$. Then
\begin{align*}
&\mathrm{MRVaR}_{(\boldsymbol{p,q})}(\mathbf{X})=\boldsymbol{\mu}+\mathbf{%
\Sigma}^{\frac{1}{2}}\frac{\boldsymbol{\delta_{\boldsymbol{p,q}}}}{F_{%
\mathbf{Y}}(\boldsymbol{\eta}_{\boldsymbol{p}},\boldsymbol{\eta}_{%
\boldsymbol{q}})},
\end{align*}
where
\begin{equation*}
\boldsymbol{\delta_{\boldsymbol{p,q}}}=\left(\delta_{1,\boldsymbol{p,q}%
},~\delta_{2,\boldsymbol{p,q}},\cdots,\delta_{n,\boldsymbol{p,q}}\right)^{T}
\end{equation*}
with
\begin{align*}
&\delta_{k,\boldsymbol{p,q}}=\frac{c_{n}}{c_{n-1,\boldsymbol{p},k}^{\ast}}F_{%
\mathbf{Y}_{\boldsymbol{p},-k}}(\boldsymbol{\eta}_{\boldsymbol{p},-k},%
\boldsymbol{\eta}_{\boldsymbol{q},-k})-\frac{c_{n}}{c_{n-1,\boldsymbol{q}%
,k}^{\ast}}F_{\mathbf{Y}_{\boldsymbol{q},-k}}(\boldsymbol{\eta}_{\boldsymbol{%
p},-k},\boldsymbol{\eta}_{\boldsymbol{q},-k}),~k=1,~2,\cdots,n,
\end{align*}
$\mathbf{Y}\sim Lo_{n}\left(\boldsymbol{0},~\boldsymbol{I_{n}}\right)$,
\begin{align*}
c_{n-1,\boldsymbol{v},k}^{\ast}&=\frac{\Gamma((n-1)/2)\exp\left\{\frac{\eta_{%
\boldsymbol{v},k}^{2}}{2}\right\}}{(2\pi)^{(n-1)/2}}\left[\int_{0}^{\infty}%
\frac{t^{(n-3)/2}\exp\{-t\}}{1+\exp\left\{-\frac{\eta_{\boldsymbol{v},k}^{2}%
}{2}\right\}\exp\{-t\}}\mathrm{d}t\right]^{-1} \\
&=\frac{\exp\left\{\frac{\eta_{\boldsymbol{v},k}^{2}}{2}\right\}}{%
(2\pi)^{(n-1)/2}\Psi_{1}^{\ast}\left(-\sqrt{2\pi}\phi(\eta_{\boldsymbol{v}%
,k}),~\frac{n-1}{2},~1\right)},
\end{align*}%
and pdf of $\mathbf{Y}_{\boldsymbol{v},-k}$:
\begin{align*}
&f_{\mathbf{Y}_{\boldsymbol{v},-k}}(\boldsymbol{t})=c_{n-1,\boldsymbol{v}%
,k}^{\ast}\frac{\exp\left\{-\frac{\boldsymbol{t}^{T}\boldsymbol{t}}{2}-\frac{%
\eta_{\boldsymbol{v},k}^{2}}{2}\right\}}{1+\exp\left\{-\frac{\boldsymbol{t}%
^{T}\boldsymbol{t}}{2}-\frac{\eta_{\boldsymbol{v},k}^{2}}{2}\right\}}, ~%
\boldsymbol{t}\in\mathbb{R}^{n-1},~\boldsymbol{v}\in\{\boldsymbol{p,q}%
\},~k=1,~2,\cdots,n.
\end{align*}
\end{corollary}

We can further simplify
\begin{equation*}
\frac{c_{n}}{c_{n-1,\boldsymbol{v},k}^{\ast}}=\frac{\Psi_{1}^{\ast}\left(-%
\sqrt{2\pi}\phi(\eta_{\boldsymbol{v},k}),~\frac{n-1}{2},~1\right)\phi(\eta_{%
\boldsymbol{v},k})}{\Psi_{2}^{\ast}\left(-1,\frac{n}{2},1\right)}%
,~k=1,2,\cdots,n,~\boldsymbol{v}\in\{\boldsymbol{p,q}\}.
\end{equation*}

\begin{corollary}
\label{co.4} Let $\mathbf{X}\sim La_{n}(\boldsymbol{\mu},~\mathbf{\Sigma})$ $%
(n\geq2)$. Then
\begin{align*}
&\mathrm{MRVaR}_{(\boldsymbol{p,q})}(\mathbf{X})=\boldsymbol{\mu}+\mathbf{%
\Sigma}^{\frac{1}{2}}\frac{\boldsymbol{\delta_{\boldsymbol{p,q}}}}{F_{%
\mathbf{Y}}(\boldsymbol{\eta}_{\boldsymbol{p}},\boldsymbol{\eta}_{%
\boldsymbol{q}})},
\end{align*}
where
\begin{equation*}
\boldsymbol{\delta_{\boldsymbol{p,q}}}=\left(\delta_{1,\boldsymbol{p,q}%
},~\delta_{2,\boldsymbol{p,q}},\cdots,\delta_{n,\boldsymbol{p,q}}\right)^{T}
\end{equation*}
with
\begin{align*}
&\delta_{k,\boldsymbol{p,q}}=\frac{c_{n}}{c_{n-1,\boldsymbol{p},k}^{\ast}}F_{%
\mathbf{Y}_{\boldsymbol{p},-k}}(\boldsymbol{\eta}_{\boldsymbol{p},-k},%
\boldsymbol{\eta}_{\boldsymbol{q},-k})-\frac{c_{n}}{c_{n-1,\boldsymbol{q}%
,k}^{\ast}}F_{\mathbf{Y}_{\boldsymbol{q},-k}}(\boldsymbol{\eta}_{\boldsymbol{%
p},-k},\boldsymbol{\eta}_{\boldsymbol{q},-k}),~k=1,~2,\cdots,n,
\end{align*}
$\mathbf{Y}\sim La_{n}\left(\boldsymbol{0},~\boldsymbol{I_{n}}\right)$,
\begin{equation*}
c_{n-1,\boldsymbol{v},k}^{\ast}=\frac{\Gamma\left(\frac{n-1}{2}\right)}{%
(2\pi)^{(n-1)/2}}\left[\int_{0}^{\infty}t^{\frac{n-3}{2}}\left(1+\sqrt{%
2t+\eta_{\boldsymbol{v},k}^{2}}\right)\exp\left\{-\sqrt{2t+\eta_{\boldsymbol{%
v},k}^{2}}\right\}\mathrm{d}t\right]^{-1},
\end{equation*}
and pdf of $\mathbf{Y}_{\boldsymbol{v},-k}$:
\begin{align*}
f_{\mathbf{Y}_{\boldsymbol{v},-k}}(\boldsymbol{t})=c_{n-1,\boldsymbol{v}%
,k}^{\ast}\left(1+\sqrt{\boldsymbol{t}^{T}\boldsymbol{t}+\eta_{\boldsymbol{v}%
,k}^{2}}\right)\exp\left\{-\sqrt{\boldsymbol{t}^{T}\boldsymbol{t}+\eta_{%
\boldsymbol{v},k}^{2}}\right\},~k=1,~2,\cdots,n,
\end{align*}
$\boldsymbol{t}\in\mathbb{R}^{n-1}$,~$\boldsymbol{v}\in\{\boldsymbol{p},%
\boldsymbol{q}\}$.
\end{corollary}

We can further simplify
\begin{equation*}
\frac{c_{n}}{c_{n-1,\boldsymbol{v},k}^{\ast}}=\frac{\Gamma\left(\frac{n}{2}%
\right)2^{(n-3)/2}}{\sqrt{\pi}\Gamma\left(\frac{n-1}{2}\right)\Gamma(n)}%
\left[\int_{0}^{\infty}t^{\frac{n-3}{2}}\left(1+\sqrt{2t+\eta_{\boldsymbol{v}%
,k}^{2}}\right)\exp\left\{-\sqrt{2t+\eta_{\boldsymbol{v},k}^{2}}\right\}%
\mathrm{d}t\right],
\end{equation*}
$~k=1,~2,\cdots,n,~\boldsymbol{v}=\boldsymbol{p},\boldsymbol{q}.$

\begin{corollary}
\label{co.5} Let $\mathbf{X}\sim PVII_{n}(\boldsymbol{\mu},~\mathbf{\Sigma}%
,~t)$ $(n\geq2)$. Then
\begin{align*}
&\mathrm{MRVaR}_{(\boldsymbol{p,q})}(\mathbf{X})=\boldsymbol{\mu}+\mathbf{%
\Sigma}^{\frac{1}{2}}\frac{\boldsymbol{\delta_{\boldsymbol{p,q}}}}{F_{%
\mathbf{Y}}(\boldsymbol{\eta}_{\boldsymbol{p}},\boldsymbol{\eta}_{%
\boldsymbol{q}})},
\end{align*}
where
\begin{equation*}
\boldsymbol{\delta_{\boldsymbol{p,q}}}=\left(\delta_{1,\boldsymbol{p,q}%
},~\delta_{2,\boldsymbol{p,q}},\cdots,\delta_{n,\boldsymbol{p,q}}\right)^{T}
\end{equation*}
with
\begin{align*}
&\delta_{k,\boldsymbol{p,q}}=\frac{c_{n}}{c_{n-1,\boldsymbol{p},k}^{\ast}}F_{%
\mathbf{Y}_{\boldsymbol{p},-k}}(\boldsymbol{\eta}_{\boldsymbol{p},-k},%
\boldsymbol{\eta}_{\boldsymbol{q},-k})-\frac{c_{n}}{c_{n-1,\boldsymbol{q}%
,k}^{\ast}}F_{\mathbf{Y}_{\boldsymbol{q},-k}}(\boldsymbol{\eta}_{\boldsymbol{%
p},-k},\boldsymbol{\eta}_{\boldsymbol{q},-k}),~k=1,~2,\cdots,n,
\end{align*}
$\mathbf{Y}\sim PVII_{n}\left(\boldsymbol{0},~\boldsymbol{I_{n}},~t\right)$,
$\mathbf{Y}_{\boldsymbol{v},-k}\sim PVII_{n-1}\left(\boldsymbol{0},~\Lambda_{%
\boldsymbol{v},k},~t-1\right),$
\begin{equation*}
\Lambda_{\boldsymbol{v},k}=\left(1+\eta_{\boldsymbol{v},k}^{2}\right)\mathbf{%
I}_{n-1}
\end{equation*}
and
\begin{equation*}
c_{n-1,\boldsymbol{v},k}^{\ast}=\frac{2\Gamma(t)}{\Gamma\left(t-\frac{n-1}{2}%
\right)\pi^{(n-1)/2}}\left(1+\eta_{\boldsymbol{v},k}^{2}\right)^{t-1},~%
\boldsymbol{v}\in\{\boldsymbol{p},\boldsymbol{q}\}.
\end{equation*}
\end{corollary}

We can further simplify
\begin{equation*}
\frac{c_{n}}{c_{n-1,\boldsymbol{v},k}^{\ast }}=\frac{\Gamma \left( t-\frac{%
n-1}{2}\right) }{2\Gamma \left( t-\frac{n}{2}\right) \sqrt{\pi }}\left(
1+\eta _{\boldsymbol{v},k}^{2}\right) ^{-(t-1)},~\boldsymbol{v}\in \{%
\boldsymbol{p},\boldsymbol{q}\}.
\end{equation*}

\subsection{RV and MRCov of elliptical distribution\label{sec:5}}

To derive the formula for MRCov, we further define the shifted cumulative
generator as
\begin{equation*}
\overline{\mathcal{G}}_{n-2,\boldsymbol{u}r,\boldsymbol{v}l}^{\ast \ast }(t)=%
\overline{\mathcal{G}}_{n}\left( t+\frac{1}{2}\eta _{\boldsymbol{u},r}^{2}+%
\frac{1}{2}\eta _{\boldsymbol{v},l}^{2}\right) ,~\boldsymbol{u},\boldsymbol{v%
}\in \{\boldsymbol{p},\boldsymbol{q}\};r\neq l;r,l\in \{1,\cdots ,n\},
\end{equation*}%
and normalizing constant
\begin{equation*}
c_{n-2,\boldsymbol{u}r,\boldsymbol{v}l}^{\ast \ast }=\frac{\Gamma (\frac{n-2%
}{2})}{(2\pi )^{(n-2)/2}}\left[ \int_{0}^{\infty }s^{(n-4)/2}\overline{%
\mathcal{G}}_{n-2,\boldsymbol{u}r,\boldsymbol{v}l}^{\ast \ast }(t)\mathrm{d}t%
\right] ^{-1}.
\end{equation*}


\begin{theorem}
\label{th.2.2} Under conditions (\ref{(v3)}), (\ref{(v5)}) and (\ref{(v6)}),
$\mathrm{RV}$ of elliptical distributed random variable $X\sim
E_{1}\left(\mu,~\sigma^2,~g_{1}\right)$ is given by
\begin{align}  \label{(v9)}
&\mathrm{RV}_{(p,q)}(X)=\frac{\sigma^{2}}{F_{Y}(\eta_{p},\eta_{q})}%
\left\{\lambda_{p,q}+\frac{c_{1}}{c_{1}^{\ast}}F_{Y^{\ast}}(\eta_{p},%
\eta_{q})-\delta_{p,q}^{2}\right\},
\end{align}
where
\begin{align*}
\lambda_{p,q}=c_{1}\left[\eta_{p}\overline{G}_{1}\left(\frac{1}{2}%
\eta_{p}^{2}\right)-\eta_{q}\overline{G}_{1}\left(\frac{1}{2}%
\eta_{q}^{2}\right)\right],
\end{align*}
$Y\sim E_{1}\left(0,~1,~g_{1}\right)$, $Y^{\ast}\sim E_{1}\left(0,~1,~%
\overline{G}_{1}\right)$, and $\delta_{p,q}$ is the same as in (\ref{(v8)}).
\end{theorem}

\begin{theorem}
\label{th.2} Under conditions (\ref{(v3)}), (\ref{(v5)}) and (\ref{(v6)}), $%
\mathrm{MRCov}$ of elliptical distributed random vector $\mathbf{X}\sim
E_{n}\left(\boldsymbol{\mu},~\mathbf{\Sigma},~g_{n}\right)$ is given by
\begin{align}  \label{(v10)}
&\mathrm{MRCov}_{(\boldsymbol{p},\boldsymbol{q})}(\mathbf{X})=\mathbf{\Sigma}%
^{\frac{1}{2}}\mathbf{\Upsilon}_{\boldsymbol{p},\boldsymbol{q}}\mathbf{\Sigma%
}^{\frac{1}{2}},
\end{align}
where
\begin{align*}
\mathbf{\Upsilon}_{\boldsymbol{p,q},i,j}=&\frac{1}{F_{\mathbf{Y}}(%
\boldsymbol{\eta_{p}},\boldsymbol{\eta_{q}})}\bigg\{\frac{c_{n}}{c_{n-2,%
\boldsymbol{p}i,\boldsymbol{p}j}^{\ast\ast}}F_{\mathbf{Y}_{\boldsymbol{p}i,%
\boldsymbol{p}j}}(\boldsymbol{\eta}_{\boldsymbol{p},-ij},\boldsymbol{\eta}_{%
\boldsymbol{q},-ij})-\frac{c_{n}}{c_{n-2,\boldsymbol{p}i,\boldsymbol{q}%
j}^{\ast\ast}}F_{\mathbf{Y}_{\boldsymbol{p}i,\boldsymbol{q}j}}(\boldsymbol{%
\eta}_{\boldsymbol{p},-ij},\boldsymbol{\eta}_{\boldsymbol{q},-ij}) \\
&+\frac{c_{n}}{c_{n-2,\boldsymbol{q}i,\boldsymbol{q}j}^{\ast\ast}}F_{\mathbf{%
Y}_{\boldsymbol{q}i,\boldsymbol{q}j}}(\boldsymbol{\eta}_{\boldsymbol{p},-ij},%
\boldsymbol{\eta}_{\boldsymbol{q},-ij}) -\frac{c_{n}}{c_{n-2,\boldsymbol{p}j,%
\boldsymbol{q}i}^{\ast\ast}}F_{\mathbf{Y}_{\boldsymbol{p}j,\boldsymbol{q}i}}(%
\boldsymbol{\eta}_{\boldsymbol{p},-ij},\boldsymbol{\eta}_{\boldsymbol{q}%
,-ij}) \bigg\} \\
&-\mathrm{MRVaR}_{(\boldsymbol{\eta_{p}},\boldsymbol{\eta_{q}})}(\mathbf{Y}%
)_{i}\mathrm{MRVaR}_{(\boldsymbol{\eta_{p}},\boldsymbol{\eta_{q}})}(\mathbf{Y%
})_{j} ,~i\neq j;~i,j\in\{1,2,\cdots,n\}
\end{align*}
and
\begin{align*}
\mathbf{\Upsilon}_{\boldsymbol{p,q},ii}=&\frac{1}{F_{\mathbf{Y}}(\boldsymbol{%
\eta}_{\boldsymbol{p}},\boldsymbol{\eta}_{\boldsymbol{q}})}\bigg\{\frac{c_{n}%
}{c_{n-1,\boldsymbol{p},i}^{\ast}}\eta_{\boldsymbol{p},i}F_{\mathbf{Y}_{%
\boldsymbol{p},-i}}(\boldsymbol{\eta}_{\boldsymbol{p},-i},\boldsymbol{\eta}_{%
\boldsymbol{q},-i}) -\frac{c_{n}}{c_{n-1,\boldsymbol{q},i}^{\ast}}\eta_{%
\boldsymbol{q},i}F_{\mathbf{Y}_{\boldsymbol{q},-i}}(\boldsymbol{\eta}_{%
\boldsymbol{p},-i},\boldsymbol{\eta}_{\boldsymbol{q},-i}) \\
&+\frac{c_{n}}{c_{n}^{\ast}}F_{\mathbf{Y}^{\ast}}(\boldsymbol{\eta}_{%
\boldsymbol{p}},\boldsymbol{\eta}_{\boldsymbol{q}})\bigg\}-[\mathrm{MRVaR}_{(%
\boldsymbol{\eta_{p}},\boldsymbol{\eta_{q}})}(\mathbf{Y})_{i}]^{2},~i=1,2,%
\cdots,n,
\end{align*}
with
\begin{align}  \label{(v11)}
&\mathrm{MRVaR}_{(\boldsymbol{\eta_{p}},\boldsymbol{\eta_{q}})}(\mathbf{Y}%
)_{k}=\frac{1}{F_{\mathbf{Y}}(\boldsymbol{\eta}_{\boldsymbol{p}},\boldsymbol{%
\eta}_{\boldsymbol{q}})}\bigg\{\frac{c_{n}}{c_{n-1,\boldsymbol{p},k}^{\ast}}%
F_{\mathbf{Y}_{\boldsymbol{p},-k}}(\boldsymbol{\eta}_{\boldsymbol{p},-k},%
\boldsymbol{\eta}_{\boldsymbol{q},-k})-\frac{c_{n}}{c_{n-1,\boldsymbol{q}%
,k}^{\ast}}F_{\mathbf{Y}_{\boldsymbol{q},-k}}(\boldsymbol{\eta}_{\boldsymbol{%
p},-k},\boldsymbol{\eta}_{\boldsymbol{q},-k})\bigg\},
\end{align}
$k=1,~2,\cdots,n,$ $\mathbf{Y}\sim E_{n}\left(\boldsymbol{0},~\boldsymbol{%
I_{n}},~g_{n}\right)$, $\mathbf{Y}^{\ast}\sim E_{n}\left(\boldsymbol{0},~%
\boldsymbol{I_{n}},~\overline{G}_{n}\right)$, $\mathbf{Y}_{\boldsymbol{u}r,%
\boldsymbol{v}l}\sim E_{n-2}(\boldsymbol{0},~\mathbf{I_{n-2}},~\overline{%
\mathcal{G}}_{n-2,\boldsymbol{u}r,\boldsymbol{v}l}^{\ast\ast})$, $%
\boldsymbol{u},\boldsymbol{v}\in\{\boldsymbol{p},\boldsymbol{q}\};r\neq
l;r,l\in\{1,\cdots,n\}$, $\mathbf{Y}_{\boldsymbol{v},-k}\sim E_{n-1}\left(%
\boldsymbol{0},~\boldsymbol{I_{n-1}},~\overline{G}_{n-1,\boldsymbol{v}%
,k}^{\ast}\right)$ .
\end{theorem}

\begin{remark}
\label{re.5} We find that $\boldsymbol{q}\rightarrow(1,\cdots,1)^{T}$ in
Theorem \ref{th.2}, one gets the formula of Theorem 2 in Landsman et al.
(2018); When $q\rightarrow 1$ in (\ref{(v9)}), it coincides with the result
of (1.7) in Furman and Landsman (2006), which is the generalization of
Theorem 2 in Jiang et al. (2016).
\end{remark}

Replacing the density generators $g_{n}(u)$, $\overline{G}_{n}(u)$ and $%
\overline{\mathcal{G}}_{n}(u)$ of Examples \ref{exa.1}-\ref{exa.5} into
Theorem \ref{th.2}, and doing some algebraic operations, we can obtain the
following Corollaries \ref{co.6}-\ref{co.10}.

\begin{corollary}
\label{co.6} Let $\mathbf{X}\sim N_{n}(\boldsymbol{\mu},~\mathbf{\Sigma})$ $%
(n\geq2)$. Then
\begin{align*}
&\mathrm{MRCov}_{(\boldsymbol{p},\boldsymbol{q})}(\mathbf{X})=\mathbf{\Sigma}%
^{\frac{1}{2}}\mathbf{\Upsilon}_{\boldsymbol{p},\boldsymbol{q}}\mathbf{\Sigma%
}^{\frac{1}{2}},
\end{align*}
where
\begin{align*}
\mathbf{\Upsilon}_{\boldsymbol{p,q},i,j}= & \frac{F_{\mathbf{Y},-ij}(%
\boldsymbol{\eta}_{\boldsymbol{p},-ij},\boldsymbol{\eta}_{\boldsymbol{q}%
,-ij})}{ F_{\mathbf{Y}}(\boldsymbol{\eta_{p}},\boldsymbol{\eta_{q}})}\bigg\{%
\phi(\eta_{\boldsymbol{p},i})\phi(\eta_{\boldsymbol{p},j})-\phi(\eta_{%
\boldsymbol{p},i})\phi(\eta_{\boldsymbol{q},j})+\phi(\eta_{\boldsymbol{q}%
,i})\phi(\eta_{\boldsymbol{q},j})-\phi(\eta_{\boldsymbol{p},j})\phi(\eta_{%
\boldsymbol{q},i})\bigg\} \\
&-\mathrm{MRVaR}_{(\boldsymbol{\eta_{p}},\boldsymbol{\eta_{q}})}(\mathbf{Y}%
)_{i}\mathrm{MRVaR}_{(\boldsymbol{\eta_{p}},\boldsymbol{\eta_{q}})}(\mathbf{Y%
})_{j} ,~i\neq j;~i,j\in\{1,2,\cdots,n\},
\end{align*}
and
\begin{align*}
\mathbf{\Upsilon}_{\boldsymbol{p,q},ii} &=\frac{F_{\mathbf{Y}_{-i}}(%
\boldsymbol{\eta}_{\boldsymbol{p},-i},\boldsymbol{\eta}_{\boldsymbol{q},-i})%
}{F_{\mathbf{Y}}(\boldsymbol{\eta}_{\boldsymbol{p}},\boldsymbol{\eta}_{%
\boldsymbol{q}})}\left[\eta_{\boldsymbol{p},i}\phi(\eta_{\boldsymbol{p}%
,i})-\eta_{\boldsymbol{q},i}\phi(\eta_{\boldsymbol{q},i})\right]+1-[\mathrm{%
MRE}_{(\boldsymbol{\eta_{p}},\boldsymbol{\eta_{q}})}(\mathbf{Y}%
)_{i}]^{2},~i=1,2,\cdots,n,
\end{align*}
with
\begin{align*}
\mathrm{MRVaR}_{(\boldsymbol{\eta_{p}},\boldsymbol{\eta_{q}})}(\mathbf{Y}%
)_{k} &=\frac{F_{\mathbf{Y}_{-k}}(\boldsymbol{\eta}_{\boldsymbol{p},-k},%
\boldsymbol{\eta}_{\boldsymbol{q},-k})}{F_{\mathbf{Y}}(\boldsymbol{\eta}_{%
\boldsymbol{p}},\boldsymbol{\eta}_{\boldsymbol{q}})}\left[\phi(\eta_{%
\boldsymbol{p},k})-\phi(\eta_{\boldsymbol{q},k})\right],
\end{align*}
$k=1,~2,\cdots,n,$ $\mathbf{Y}\sim N_{n}\left(\boldsymbol{0},~\boldsymbol{%
I_{n}}\right)$, $\mathbf{Y}_{-k}\sim N_{n-1}(\boldsymbol{0},~\mathbf{I_{n-1}}%
)$, $\mathbf{Y}_{-ij}\sim N_{n-2}(\boldsymbol{0},~\mathbf{I_{n-2}})$.
\end{corollary}

\begin{corollary}
\label{co.7} Let $\mathbf{X}\sim St_{n}(\boldsymbol{\mu},~\mathbf{\Sigma},m)
$ $(n\geq2)$. Then
\begin{align*}
&\mathrm{MRCov}_{(\boldsymbol{p},\boldsymbol{q})}(\mathbf{X})=\mathbf{\Sigma}%
^{\frac{1}{2}}\mathbf{\Upsilon}_{\boldsymbol{p},\boldsymbol{q}}\mathbf{\Sigma%
}^{\frac{1}{2}},
\end{align*}
where
\begin{align*}
\mathbf{\Upsilon}_{\boldsymbol{p,q},i,j}=&\frac{1}{F_{\mathbf{Y}}(%
\boldsymbol{\eta_{p}},\boldsymbol{\eta_{q}})}\bigg\{\frac{c_{n}}{c_{n-2,%
\boldsymbol{p}i,\boldsymbol{p}j}^{\ast\ast}}F_{\mathbf{Y}_{\boldsymbol{p}i,%
\boldsymbol{p}j}}(\boldsymbol{\eta}_{\boldsymbol{p},-ij},\boldsymbol{\eta}_{%
\boldsymbol{q},-ij})-\frac{c_{n}}{c_{n-2,\boldsymbol{p}i,\boldsymbol{q}%
j}^{\ast\ast}}F_{\mathbf{Y}_{\boldsymbol{p}i,\boldsymbol{q}j}}(\boldsymbol{%
\eta}_{\boldsymbol{p},-ij},\boldsymbol{\eta}_{\boldsymbol{q},-ij}) \\
&+\frac{c_{n}}{c_{n-2,\boldsymbol{q}i,\boldsymbol{q}j}^{\ast\ast}}F_{\mathbf{%
Y}_{\boldsymbol{q}i,\boldsymbol{q}j}}(\boldsymbol{\eta}_{\boldsymbol{p},-ij},%
\boldsymbol{\eta}_{\boldsymbol{q},-ij}) -\frac{c_{n}}{c_{n-2,\boldsymbol{p}j,%
\boldsymbol{q}i}^{\ast\ast}}F_{\mathbf{Y}_{\boldsymbol{p}j,\boldsymbol{q}i}}(%
\boldsymbol{\eta}_{\boldsymbol{p},-ij},\boldsymbol{\eta}_{\boldsymbol{q}%
,-ij}) \bigg\} \\
&-\mathrm{MRVaR}_{(\boldsymbol{\eta_{p}},\boldsymbol{\eta_{q}})}(\mathbf{Y}%
)_{i}\mathrm{MRVaR}_{(\boldsymbol{\eta_{p}},\boldsymbol{\eta_{q}})}(\mathbf{Y%
})_{j} ,~i\neq j,~i,j\in\{1,2,\cdots,n\},
\end{align*}
and
\begin{align*}
\mathbf{\Upsilon}_{\boldsymbol{p,q},ii}&=\frac{1}{F_{\mathbf{Y}}(\boldsymbol{%
\eta}_{\boldsymbol{p}},\boldsymbol{\eta}_{\boldsymbol{q}})}\bigg\{\frac{c_{n}%
}{c_{n-1,\boldsymbol{p},i}^{\ast}}\eta_{\boldsymbol{p},i}F_{\mathbf{Y}_{%
\boldsymbol{p},-i}}(\boldsymbol{\eta}_{\boldsymbol{p},-i},\boldsymbol{\eta}_{%
\boldsymbol{q},-i}) -\frac{c_{n}}{c_{n-1,\boldsymbol{q},i}^{\ast}}\eta_{%
\boldsymbol{q},i}F_{\mathbf{Y}_{\boldsymbol{q},-i}}(\boldsymbol{\eta}_{%
\boldsymbol{p},-i},\boldsymbol{\eta}_{\boldsymbol{q},-i}) \\
&~~~+\frac{c_{n}}{c_{n}^{\ast}}F_{\mathbf{Y}^{\ast}}(\boldsymbol{\eta}_{%
\boldsymbol{p}},\boldsymbol{\eta}_{\boldsymbol{q}})\bigg\}-[\mathrm{MRVaR}_{(%
\boldsymbol{\eta_{p}},\boldsymbol{\eta_{q}})}(\mathbf{Y})_{i}]^{2},~i=1,2,%
\cdots,n,
\end{align*}
with
\begin{align*}
&\mathrm{MRVaR}_{(\boldsymbol{\eta_{p}},\boldsymbol{\eta_{q}})}(\mathbf{Y}%
)_{k}=\frac{1}{F_{\mathbf{Y}}(\boldsymbol{\eta}_{\boldsymbol{p}},\boldsymbol{%
\eta}_{\boldsymbol{q}})}\bigg\{\frac{c_{n}}{c_{n-1,\boldsymbol{p},k}^{\ast}}%
F_{\mathbf{Y}_{\boldsymbol{p},-k}}(\boldsymbol{\eta}_{\boldsymbol{p},-k},%
\boldsymbol{\eta}_{\boldsymbol{q},-k})-\frac{c_{n}}{c_{n-1,\boldsymbol{q}%
,k}^{\ast}}F_{\mathbf{Y}_{\boldsymbol{q},-k}}(\boldsymbol{\eta}_{\boldsymbol{%
p},-k},\boldsymbol{\eta}_{\boldsymbol{q},-k})\bigg\},
\end{align*}
$k=1,~2,\cdots,n,$ $\mathbf{Y}\sim St_{n}\left(\boldsymbol{0},~\boldsymbol{%
I_{n}},~m\right)$, $\mathbf{Y}^{\ast}\sim E_{n}\left(\boldsymbol{0},~%
\boldsymbol{I_{n}},~\overline{G}_{n}\right)$, $\mathbf{Y}_{\boldsymbol{u}r,%
\boldsymbol{v}l}\sim St_{n-2}\left(\boldsymbol{0},~\Delta_{n-2,\boldsymbol{u}%
r,\boldsymbol{v}l},~m-2\right)$, \newline
$\mathbf{Y}_{\boldsymbol{v},-k}\sim St_{n-1}\left(\boldsymbol{0},~\Delta_{%
\boldsymbol{v},k},~m-1\right),$
\begin{equation*}
\Delta_{n-2,\boldsymbol{u}r,\boldsymbol{v}l}=\left(\frac{m+\eta_{\boldsymbol{%
u}r}^{2}+\eta_{\boldsymbol{v}l}^{2}}{m-2}\right)\mathbf{I}_{n-2},
\end{equation*}
\begin{equation*}
c_{n-2,\boldsymbol{u}r,\boldsymbol{v}l}^{\ast\ast}=\frac{\Gamma\left(\frac{%
m+n-4}{2}\right)(m+n-2)(m+n-4)}{\Gamma\left(\frac{m-2}{2}\right)\pi^{\frac{%
n-2}{2}}m^{\frac{m+n}{2}}}\left(m+\eta_{\boldsymbol{u}r}^{2}+\eta_{%
\boldsymbol{v}l}^{2}\right)^{\frac{m+n-4}{2}},
\end{equation*}
$\boldsymbol{u},\boldsymbol{v}\in\{\boldsymbol{p},\boldsymbol{q}\};r\neq
l;r,l\in\{1,\cdots,n\},$ $c_{n-1,\boldsymbol{u},k}^{\ast}$ and $\Delta_{%
\boldsymbol{v},k}$ are the same as those in Corollary \ref{co.2}.
\end{corollary}

We can further simplify
\begin{equation*}
\frac{c_{n}}{c_{n-2,\boldsymbol{u}r,\boldsymbol{v}l}^{\ast\ast}}=\frac{%
(m+n)m^{\frac{m-2}{2}}}{2(m+n-4)\pi}\left(m+\eta_{\boldsymbol{u}r}^{2}+\eta_{%
\boldsymbol{v}l}^{2}\right)^{-(m+n-4)/2}
\end{equation*}
and
\begin{equation*}
\frac{c_{n}}{c_{n}^{\ast}} =\frac{p\Gamma(\frac{m+n}{2})B(\frac{n}{2},~\frac{%
m-2}{2})}{(m+n-2)\Gamma(\frac{m}{2})\Gamma(\frac{n}{2})}=\frac{m}{m-2},~%
\mathrm{if}~m>2.
\end{equation*}

\begin{corollary}
\label{co.8} Let $\mathbf{X}\sim Lo_{n}(\boldsymbol{\mu},~\mathbf{\Sigma})$ $%
(n\geq2)$. Then
\begin{align*}
&\mathrm{MRCov}_{(\boldsymbol{p},\boldsymbol{q})}(\mathbf{X})=\mathbf{\Sigma}%
^{\frac{1}{2}}\mathbf{\Upsilon}_{\boldsymbol{p},\boldsymbol{q}}\mathbf{\Sigma%
}^{\frac{1}{2}},
\end{align*}
where
\begin{align*}
\mathbf{\Upsilon}_{\boldsymbol{p,q},i,j}=& \frac{1}{F_{\mathbf{Y}}(%
\boldsymbol{\eta_{p}},\boldsymbol{\eta_{q}})}\bigg\{\frac{c_{n}}{c_{n-2,%
\boldsymbol{p}i,\boldsymbol{p}j}^{\ast\ast}}F_{\mathbf{Y}_{\boldsymbol{p}i,%
\boldsymbol{p}j}}(\boldsymbol{\eta}_{\boldsymbol{p},-ij},\boldsymbol{\eta}_{%
\boldsymbol{q},-ij})-\frac{c_{n}}{c_{n-2,\boldsymbol{p}i,\boldsymbol{q}%
j}^{\ast\ast}}F_{\mathbf{Y}_{\boldsymbol{p}i,\boldsymbol{q}j}}(\boldsymbol{%
\eta}_{\boldsymbol{p},-ij},\boldsymbol{\eta}_{\boldsymbol{q},-ij}) \\
&+\frac{c_{n}}{c_{n-2,\boldsymbol{q}i,\boldsymbol{q}j}^{\ast\ast}}F_{\mathbf{%
Y}_{\boldsymbol{q}i,\boldsymbol{q}j}}(\boldsymbol{\eta}_{\boldsymbol{p},-ij},%
\boldsymbol{\eta}_{\boldsymbol{q},-ij}) -\frac{c_{n}}{c_{n-2,\boldsymbol{p}j,%
\boldsymbol{q}i}^{\ast\ast}}F_{\mathbf{Y}_{\boldsymbol{p}j,\boldsymbol{q}i}}(%
\boldsymbol{\eta}_{\boldsymbol{p},-ij},\boldsymbol{\eta}_{\boldsymbol{q}%
,-ij}) \bigg\} \\
&-\mathrm{MRVaR}_{(\boldsymbol{\eta_{p}},\boldsymbol{\eta_{q}})}(\mathbf{Y}%
)_{i}\mathrm{MRVaR}_{(\boldsymbol{\eta_{p}},\boldsymbol{\eta_{q}})}(\mathbf{Y%
})_{j} ,~i\neq j;~i,j\in\{1,2,\cdots,n\},
\end{align*}
and
\begin{align*}
\mathbf{\Upsilon}_{\boldsymbol{p,q},ii}&=\frac{1}{F_{\mathbf{Y}}(\boldsymbol{%
\eta}_{\boldsymbol{p}},\boldsymbol{\eta}_{\boldsymbol{q}})}\bigg\{\frac{c_{n}%
}{c_{n-1,\boldsymbol{p},i}^{\ast}}\eta_{\boldsymbol{p},i}F_{\mathbf{Y}_{%
\boldsymbol{p},-i}}(\boldsymbol{\eta}_{\boldsymbol{p},-i},\boldsymbol{\eta}_{%
\boldsymbol{q},-i}) -\frac{c_{n}}{c_{n-1,\boldsymbol{q},i}^{\ast}}\eta_{%
\boldsymbol{q},i}F_{\mathbf{Y}_{\boldsymbol{q},-i}}(\boldsymbol{\eta}_{%
\boldsymbol{p},-i},\boldsymbol{\eta}_{\boldsymbol{q},-i}) \\
&~~~+\frac{\Psi_{1}^{\ast}(-1,\frac{n}{2},1)}{\Psi_{2}^{\ast}(-1,\frac{n}{2}%
,1)}F_{\mathbf{Y}^{\ast}}(\boldsymbol{\eta}_{\boldsymbol{p}},\boldsymbol{\eta%
}_{\boldsymbol{q}})\bigg\}-[\mathrm{MRVaR}_{(\boldsymbol{\eta_{p}},%
\boldsymbol{\eta_{q}})}(\mathbf{Y})_{i}]^{2},~i=1,2,\cdots,n,
\end{align*}
with
\begin{align*}
&\mathrm{MRVaR}_{(\boldsymbol{\eta_{p}},\boldsymbol{\eta_{q}})}(\mathbf{Y}%
)_{k}=\frac{1}{F_{\mathbf{Y}}(\boldsymbol{\eta}_{\boldsymbol{p}},\boldsymbol{%
\eta}_{\boldsymbol{q}})}\bigg\{\frac{c_{n}}{c_{n-1,\boldsymbol{p},k}^{\ast}}%
F_{\mathbf{Y}_{\boldsymbol{p},-k}}(\boldsymbol{\eta}_{\boldsymbol{p},-k},%
\boldsymbol{\eta}_{\boldsymbol{q},-k})-\frac{c_{n}}{c_{n-1,\boldsymbol{q}%
,k}^{\ast}}F_{\mathbf{Y}_{\boldsymbol{q},-k}}(\boldsymbol{\eta}_{\boldsymbol{%
p},-k},\boldsymbol{\eta}_{\boldsymbol{q},-k})\bigg\},
\end{align*}
$k=1,~2,\cdots,n,$ $\mathbf{Y}\sim Lo_{n}\left(\boldsymbol{0},~\boldsymbol{%
I_{n}}\right)$, $\mathbf{Y}^{\ast}\sim E_{n}\left(\boldsymbol{0},~%
\boldsymbol{I_{n}},~\overline{G}_{n}\right)$, the pdf of $\mathbf{Y}_{%
\boldsymbol{u}r,\boldsymbol{v}l}$:
\begin{equation*}
f_{\mathbf{Y}_{\boldsymbol{u}r,\boldsymbol{v}l}}(\boldsymbol{t})=c_{n-2,%
\boldsymbol{u}r,\boldsymbol{v}l}^{\ast\ast}\ln\left[1+\exp\left(-\frac{1}{2}%
\boldsymbol{t}^{T}\boldsymbol{t}-\frac{1}{2}\eta_{\boldsymbol{u},r}^{2}-%
\frac{1}{2}\eta_{\boldsymbol{v},l}^{2}\right)\right],~\boldsymbol{t}\in
\mathbb{R}^{n-2},
\end{equation*}
the normalizing constant
\begin{align*}
&c_{n-2,\boldsymbol{u}r,\boldsymbol{v}l}^{\ast\ast}=\frac{\Gamma((n-2)/2)}{%
(2\pi)^{(n-2)/2}}\left\{\int_{0}^{\infty}t^{(n-4)/2}\ln\left[1+\exp\left(-%
\frac{1}{2}\eta_{\boldsymbol{u},r}^{2}-\frac{1}{2}\eta_{\boldsymbol{v}%
,l}^{2}\right)\exp(-t)\right]\mathrm{d}t\right\}^{-1},
\end{align*}
$\boldsymbol{u},\boldsymbol{v}\in\{\boldsymbol{p},\boldsymbol{q}\};r\neq
l;r,l\in\{1,\cdots,n\},$ and $\mathbf{Y}_{\boldsymbol{v},-k}$ and $c_{n-1,%
\boldsymbol{p},k}^{\ast}$ are the same as those in Corollary \ref{co.3}.
\end{corollary}

We can further simplify
\begin{equation*}
\frac{c_{n}}{c_{n-2,\boldsymbol{u}r,\boldsymbol{v}l}^{\ast\ast}}=\frac{%
\int_{0}^{\infty}t^{(n-4)/2}\ln\left[1+\exp\left(-\frac{1}{2}\eta_{%
\boldsymbol{u},r}^{2}-\frac{1}{2}\eta_{\boldsymbol{v},l}^{2}\right)\exp(-t)%
\right]\mathrm{d}t}{2\Gamma((n-1)/2)\pi\Psi_{2}^{\ast}(-1,\frac{n}{2},1)}.
\end{equation*}

\begin{corollary}
\label{co.9} Let $\mathbf{X}\sim La_{n}(\boldsymbol{\mu},~\mathbf{\Sigma})$ $%
(n\geq2)$. Then
\begin{align*}
&\mathrm{MRCov}_{(\boldsymbol{p},\boldsymbol{q})}(\mathbf{X})=\mathbf{\Sigma}%
^{\frac{1}{2}}\mathbf{\Upsilon}_{\boldsymbol{p},\boldsymbol{q}}\mathbf{\Sigma%
}^{\frac{1}{2}},
\end{align*}
where
\begin{align*}
\mathbf{\Upsilon}_{\boldsymbol{p,q},i,j}=&\frac{1}{F_{\mathbf{Y}}(%
\boldsymbol{\eta_{p}},\boldsymbol{\eta_{q}})}\bigg\{\frac{c_{n}}{c_{n-2,%
\boldsymbol{p}i,\boldsymbol{p}j}^{\ast\ast}}F_{\mathbf{Y}_{\boldsymbol{p}i,%
\boldsymbol{p}j}}(\boldsymbol{\eta}_{\boldsymbol{p},-ij},\boldsymbol{\eta}_{%
\boldsymbol{q},-ij})-\frac{c_{n}}{c_{n-2,\boldsymbol{p}i,\boldsymbol{q}%
j}^{\ast\ast}}F_{\mathbf{Y}_{\boldsymbol{p}i,\boldsymbol{q}j}}(\boldsymbol{%
\eta}_{\boldsymbol{p},-ij},\boldsymbol{\eta}_{\boldsymbol{q},-ij}) \\
&+\frac{c_{n}}{c_{n-2,\boldsymbol{q}i,\boldsymbol{q}j}^{\ast\ast}}F_{\mathbf{%
Y}_{\boldsymbol{q}i,\boldsymbol{q}j}}(\boldsymbol{\eta}_{\boldsymbol{p},-ij},%
\boldsymbol{\eta}_{\boldsymbol{q},-ij}) -\frac{c_{n}}{c_{n-2,\boldsymbol{p}j,%
\boldsymbol{q}i}^{\ast\ast}}F_{\mathbf{Y}_{\boldsymbol{p}j,\boldsymbol{q}i}}(%
\boldsymbol{\eta}_{\boldsymbol{p},-ij},\boldsymbol{\eta}_{\boldsymbol{q}%
,-ij}) \bigg\} \\
&-\mathrm{MRVaR}_{(\boldsymbol{\eta_{p}},\boldsymbol{\eta_{q}})}(\mathbf{Y}%
)_{i}\mathrm{MRVaR}_{(\boldsymbol{\eta_{p}},\boldsymbol{\eta_{q}})}(\mathbf{Y%
})_{j} ,~i\neq j,~i,j\in\{1,2,\cdots,n\},
\end{align*}
and
\begin{align*}
\mathbf{\Upsilon}_{\boldsymbol{p,q},ii}&=\frac{1}{F_{\mathbf{Y}}(\boldsymbol{%
\eta}_{\boldsymbol{p}},\boldsymbol{\eta}_{\boldsymbol{q}})}\bigg\{\frac{c_{n}%
}{c_{n-1,\boldsymbol{p},i}^{\ast}}\eta_{\boldsymbol{p},i}F_{\mathbf{Y}_{%
\boldsymbol{p},-i}}(\boldsymbol{\eta}_{\boldsymbol{p},-i},\boldsymbol{\eta}_{%
\boldsymbol{q},-i}) -\frac{c_{n}}{c_{n-1,\boldsymbol{q},i}^{\ast}}\eta_{%
\boldsymbol{q},i}F_{\mathbf{Y}_{\boldsymbol{q},-i}}(\boldsymbol{\eta}_{%
\boldsymbol{p},-i},\boldsymbol{\eta}_{\boldsymbol{q},-i}) \\
&~~~+(n+1)F_{\mathbf{Y}^{\ast}}(\boldsymbol{\eta}_{\boldsymbol{p}},%
\boldsymbol{\eta}_{\boldsymbol{q}})\bigg\}-[\mathrm{MRVaR}_{(\boldsymbol{%
\eta_{p}},\boldsymbol{\eta_{q}})}(\mathbf{Y})_{i}]^{2},~i=1,2,\cdots,n,
\end{align*}
with
\begin{align*}
&\mathrm{MRVaR}_{(\boldsymbol{\eta_{p}},\boldsymbol{\eta_{q}})}(\mathbf{Y}%
)_{k}=\frac{1}{F_{\mathbf{Y}}(\boldsymbol{\eta}_{\boldsymbol{p}},\boldsymbol{%
\eta}_{\boldsymbol{q}})}\bigg\{\frac{c_{n}}{c_{n-1,\boldsymbol{p},k}^{\ast}}%
F_{\mathbf{Y}_{\boldsymbol{p},-k}}(\boldsymbol{\eta}_{\boldsymbol{p},-k},%
\boldsymbol{\eta}_{\boldsymbol{q},-k})-\frac{c_{n}}{c_{n-1,\boldsymbol{q}%
,k}^{\ast}}F_{\mathbf{Y}_{\boldsymbol{q},-k}}(\boldsymbol{\eta}_{\boldsymbol{%
p},-k},\boldsymbol{\eta}_{\boldsymbol{q},-k})\bigg\},
\end{align*}
$k=1,~2,\cdots,n,$ $\mathbf{Y}\sim La_{n}\left(\boldsymbol{0},~\boldsymbol{%
I_{n}}\right)$, $\mathbf{Y}^{\ast}\sim E_{n}\left(\boldsymbol{0},~%
\boldsymbol{I_{n}},~\overline{G}_{n}\right)$, the pdf of $\mathbf{Y}_{%
\boldsymbol{u}r,\boldsymbol{v}l}$:
\begin{align*}
f_{\mathbf{Y}_{\boldsymbol{u}r,\boldsymbol{v}l}}(\boldsymbol{t})=&c_{n-2,%
\boldsymbol{u}r,\boldsymbol{v}l}^{\ast\ast}\left[\left(1+\frac{3}{\sqrt{2}}%
\right)\left(\boldsymbol{t}^{T}\boldsymbol{t}+\eta_{\boldsymbol{u}%
,r}^{2}+\eta_{\boldsymbol{v},l}^{2}\right)+3\right]\exp(-\sqrt{\boldsymbol{t}%
^{T}\boldsymbol{t}+\eta_{\boldsymbol{u},r}^{2}+\eta_{\boldsymbol{v},l}^{2}}%
),~\boldsymbol{t}\in \mathbb{R}^{n-2},
\end{align*}
and the normalizing constant
\begin{align*}
c_{n-2,\boldsymbol{u}r,\boldsymbol{v}l}^{\ast\ast}=&\frac{\Gamma\left(\frac{%
n-2}{2}\right)}{(2\pi)^{(n-2)/2}}\bigg\{\int_{0}^{\infty}t^{\frac{n-4}{2}}%
\left[(2+3\sqrt{2})\left(t+\frac{\eta_{\boldsymbol{u},r}^{2}}{2}+\frac{\eta_{%
\boldsymbol{v},l}^{2}}{2}\right)+3\right]\exp\left(-\sqrt{2t+\eta_{%
\boldsymbol{u},r}^{2}+\eta_{\boldsymbol{v},l}^{2}}\right)\mathrm{d}t\bigg\}%
^{-1}.
\end{align*}
$\boldsymbol{u},\boldsymbol{v}\in\{\boldsymbol{p},\boldsymbol{q}\};r\neq
l;r,l\in\{1,\cdots,n\},$ and $c_{n-1,\boldsymbol{v},k}^{\ast}$ and $\mathbf{Y%
}_{\boldsymbol{v},-k}$ are the same as those in Corollary \ref{co.4}.
\end{corollary}

We can further simplify
\begin{align*}
\frac{c_{n}}{c_{n-2,\boldsymbol{u}r,\boldsymbol{v}l}^{\ast\ast}}=&\frac{%
2^{(n-6)/2}}{\Gamma(n-1)\pi}\bigg\{\int_{0}^{\infty}t^{\frac{n-4}{2}}\left[%
(2+3\sqrt{2})\left(t+\frac{\eta_{\boldsymbol{u},r}^{2}}{2}+\frac{\eta_{%
\boldsymbol{v},l}^{2}}{2}\right)+3\right]\exp\left(-\sqrt{2t+\eta_{%
\boldsymbol{u},r}^{2}+\eta_{\boldsymbol{v},l}^{2}}\right)\mathrm{d}t\bigg\}.
\end{align*}

\begin{corollary}
\label{co.10} Let $\mathbf{X}\sim PVII_{n}(\boldsymbol{\mu},~\mathbf{\Sigma}%
,t)$ $(n\geq2)$. Then
\begin{align*}
&\mathrm{MRCov}_{(\boldsymbol{p},\boldsymbol{q})}(\mathbf{X})=\mathbf{\Sigma}%
^{\frac{1}{2}}\mathbf{\Upsilon}_{\boldsymbol{p},\boldsymbol{q}}\mathbf{\Sigma%
}^{\frac{1}{2}},
\end{align*}
where
\begin{align*}
\mathbf{\Upsilon}_{\boldsymbol{p,q},i,j}=&\frac{1}{F_{\mathbf{Y}}(%
\boldsymbol{\eta_{p}},\boldsymbol{\eta_{q}})}\bigg\{\frac{c_{n}}{c_{n-2,%
\boldsymbol{p}i,\boldsymbol{p}j}^{\ast\ast}}F_{\mathbf{Y}_{\boldsymbol{p}i,%
\boldsymbol{p}j}}(\boldsymbol{\eta}_{\boldsymbol{p},-ij},\boldsymbol{\eta}_{%
\boldsymbol{q},-ij})-\frac{c_{n}}{c_{n-2,\boldsymbol{p}i,\boldsymbol{q}%
j}^{\ast\ast}}F_{\mathbf{Y}_{\boldsymbol{p}i,\boldsymbol{q}j}}(\boldsymbol{%
\eta}_{\boldsymbol{p},-ij},\boldsymbol{\eta}_{\boldsymbol{q},-ij}) \\
&+\frac{c_{n}}{c_{n-2,\boldsymbol{q}i,\boldsymbol{q}j}^{\ast\ast}}F_{\mathbf{%
Y}_{\boldsymbol{q}i,\boldsymbol{q}j}}(\boldsymbol{\eta}_{\boldsymbol{p},-ij},%
\boldsymbol{\eta}_{\boldsymbol{q},-ij}) -\frac{c_{n}}{c_{n-2,\boldsymbol{p}j,%
\boldsymbol{q}i}^{\ast\ast}}F_{\mathbf{Y}_{\boldsymbol{p}j,\boldsymbol{q}i}}(%
\boldsymbol{\eta}_{\boldsymbol{p},-ij},\boldsymbol{\eta}_{\boldsymbol{q}%
,-ij}) \bigg\} \\
&-\mathrm{MRVaR}_{(\boldsymbol{\eta_{p}},\boldsymbol{\eta_{q}})}(\mathbf{Y}%
)_{i}\mathrm{MRVaR}_{(\boldsymbol{\eta_{p}},\boldsymbol{\eta_{q}})}(\mathbf{Y%
})_{j} ,~i\neq j;~i,j\in\{1,2,\cdots,n\},
\end{align*}
and
\begin{align*}
\mathbf{\Upsilon}_{\boldsymbol{p,q},ii}&=\frac{1}{F_{\mathbf{Y}}(\boldsymbol{%
\eta}_{\boldsymbol{p}},\boldsymbol{\eta}_{\boldsymbol{q}})}\bigg\{\frac{c_{n}%
}{c_{n-1,\boldsymbol{p},i}^{\ast}}\eta_{\boldsymbol{p},i}F_{\mathbf{Y}_{%
\boldsymbol{p},-i}}(\boldsymbol{\eta}_{\boldsymbol{p},-i},\boldsymbol{\eta}_{%
\boldsymbol{q},-i}) -\frac{c_{n}}{c_{n-1,\boldsymbol{q},i}^{\ast}}\eta_{%
\boldsymbol{q},i}F_{\mathbf{Y}_{\boldsymbol{q},-i}}(\boldsymbol{\eta}_{%
\boldsymbol{p},-i},\boldsymbol{\eta}_{\boldsymbol{q},-i}) \\
&~~~+\frac{c_{n}}{c_{n}^{\ast}}F_{\mathbf{Y}^{\ast}}(\boldsymbol{\eta}_{%
\boldsymbol{p}},\boldsymbol{\eta}_{\boldsymbol{q}})\bigg\}-[\mathrm{MRVaR}_{(%
\boldsymbol{\eta_{p}},\boldsymbol{\eta_{q}})}(\mathbf{Y})_{i}]^{2},~i=1,2,%
\cdots,n,
\end{align*}
with
\begin{align*}
&\mathrm{MRVaR}_{(\boldsymbol{\eta_{p}},\boldsymbol{\eta_{q}})}(\mathbf{Y}%
)_{k}=\frac{1}{F_{\mathbf{Y}}(\boldsymbol{\eta}_{\boldsymbol{p}},\boldsymbol{%
\eta}_{\boldsymbol{q}})}\bigg\{\frac{c_{n}}{c_{n-1,\boldsymbol{p},k}^{\ast}}%
F_{\mathbf{Y}_{\boldsymbol{p},-k}}(\boldsymbol{\eta}_{\boldsymbol{p},-k},%
\boldsymbol{\eta}_{\boldsymbol{q},-k})-\frac{c_{n}}{c_{n-1,\boldsymbol{q}%
,k}^{\ast}}F_{\mathbf{Y}_{\boldsymbol{q},-k}}(\boldsymbol{\eta}_{\boldsymbol{%
p},-k},\boldsymbol{\eta}_{\boldsymbol{q},-k})\bigg\},
\end{align*}
$k=1,~2,\cdots,n,$ $\mathbf{Y}\sim PVII_{n}\left(\boldsymbol{0},~\boldsymbol{%
I_{n}},~m\right)$, $\mathbf{Y}^{\ast}\sim E_{n}\left(\boldsymbol{0},~%
\boldsymbol{I_{n}},~\overline{G}_{n}\right)$, $\mathbf{Y}_{\boldsymbol{u}r,%
\boldsymbol{v}l}\sim PVII_{n-2}\left(\boldsymbol{0},~\Lambda_{n-2,%
\boldsymbol{u}r,\boldsymbol{v}l},~t-2\right),$ \newline
$\mathbf{Y}_{\boldsymbol{v},-k}\sim PVII_{n-1}\left(\boldsymbol{0},~\Lambda_{%
\boldsymbol{v},k},~t-1\right),$ $\Lambda_{n-2,\boldsymbol{u}r,\boldsymbol{v}%
l}=\left(1+\eta_{\boldsymbol{u}r}^{2}+\eta_{\boldsymbol{v}l}^{2}\right)%
\mathbf{I}_{n-2},$
\begin{equation*}
c_{n-2,\boldsymbol{u}r,\boldsymbol{v}l}^{\ast\ast}=\frac{4\Gamma\left(t%
\right)}{\Gamma\left(t-\frac{n+2}{2}\right)\pi^{\frac{n-2}{2}}}\left(1+\eta_{%
\boldsymbol{u}r}^{2}+\eta_{\boldsymbol{v}l}^{2}\right)^{t-2},
\end{equation*}
$\boldsymbol{u},\boldsymbol{v}\in\{\boldsymbol{p},\boldsymbol{q}\};r\neq
l;r,l\in\{1,\cdots,n\},$ and $\mathbf{Y}_{\boldsymbol{v},-k}$, $c_{n-1,%
\boldsymbol{v},k}^{\ast}$ and $\Lambda_{\boldsymbol{v},k}$ are the same as
those in Corollary \ref{co.5}.
\end{corollary}

We can further simplify
\begin{equation*}
\frac{c_{n}}{c_{n-2,\boldsymbol{u}r,\boldsymbol{v}l}^{\ast\ast}}=\frac{1}{%
4(t-(n+2)/2)\pi}\left(1+\eta_{\boldsymbol{u}r}^{2}+\eta_{\boldsymbol{v}%
l}^{2}\right)^{-(t-2)}
\end{equation*}
and
\begin{equation*}
\frac{c_{n}}{c_{n}^{\ast}} =\frac{1}{2t-n-2},~\mathrm{if}~t>n/2+1.
\end{equation*}

\subsection{Extensions to the class of log-elliptical distributions\label%
{sec:7aa}}

The family of log-elliptical (LE) distributions is extensively used in
quantitative finance to model financial returns and losses (see Valdez et
al., 2009; Landsman et al., 2016b; Landsman and Shushi, 2021). Let $\mathbf{Z%
} \sim LE_{n}(\boldsymbol{\mu},\mathbf{\Sigma},g_{n})$, we note that $%
\mathrm{ln}\mathbf{Z}\sim E_{n}(\boldsymbol{\mu},\mathbf{\Sigma},g_{n}).$
Letting $\mathbf{P}=\mathbf{\Sigma}^{-1/2}(\mathrm{ln}\mathbf{Z}-\boldsymbol{%
\mu})$. For the ease of presentation, random vectors $\mathbf{P}_{k}\dagger$
and $\mathbf{P}_{k,j}\ddagger$ with the pdfs
\begin{align}  \label{(v12)}
f_{\mathbf{P}_{k}\dagger}(\boldsymbol{x})=\psi\left(-\frac{1}{2}%
\sigma_{kk}\right)^{-1}\exp\left\{\mathbf{\Sigma}_{k}^{1/2}\boldsymbol{x}%
\right\}c_{n}g_{n}\left(\frac{1}{2}\boldsymbol{x}^{T}\boldsymbol{x}\right)
\end{align}
and
\begin{align}  \label{(v13)}
f_{\mathbf{P}_{k,j}\ddagger}(\boldsymbol{x})=\psi\left(-\frac{1}{2}%
(\sigma_{kk}+2\sigma_{kj}+\sigma_{jj})\right)^{-1}\exp\left\{\left(\mathbf{%
\Sigma}_{k}^{1/2}+\mathbf{\Sigma}_{j}^{1/2}\right)\boldsymbol{x}%
\right\}c_{n}g_{n}\left(\frac{1}{2}\boldsymbol{x}^{T}\boldsymbol{x}\right)
\end{align}
are introduced, respectively (see Landsman and Shushi, 2021).

\begin{proposition}
\label{th.5} Let $\mathbf{Z}\sim LE_{n}(\boldsymbol{\mu },\mathbf{\Sigma }%
,g_{n})$ be a $n\times 1$ random vector of risks with characteristic
generator $\psi $, we have
\begin{equation}  \label{(v14)}
(\mathrm{i})~\mathrm{MRVaR}_{(\boldsymbol{p,q})}(\mathbf{Z})=\mathrm{e}^{%
\boldsymbol{\mu }}\circ \boldsymbol{\psi }\circ \boldsymbol{\delta _{%
\boldsymbol{p,q}}},
\end{equation}%
\begin{equation}  \label{(v15)}
(\mathrm{ii})~\mathrm{MRCov}_{(\boldsymbol{p},\boldsymbol{q})}(\mathbf{Z}%
)=(d_{kj})_{k,j=1}^{n},
\end{equation}%
where
\begin{equation*}
d_{kj}=\mathrm{e}^{\mu _{k}+\mu _{j}}\left[ \frac{\psi_{n} \left( -\frac{1}{2%
}(\sigma _{kk}+2\sigma _{kj}+\sigma _{jj})\right) F_{\mathbf{P}%
_{k,j}\ddagger }(\boldsymbol{\zeta _{p}},\boldsymbol{\zeta _{q}})}{F_{%
\mathbf{Z}}(\mathrm{VaR}_{\boldsymbol{p}}(\mathbf{Z}),\mathrm{VaR}_{%
\boldsymbol{q}}(\mathbf{Z}))}-\frac{\psi_{n} \left( -\frac{1}{2}\sigma
_{kk}\right) \psi_{n} \left( -\frac{1}{2}\sigma _{jj}\right) F_{\mathbf{P}%
_{k}\dagger }(\boldsymbol{\zeta _{p}},\boldsymbol{\zeta _{q}})F_{\mathbf{P}%
_{j}\dagger }(\boldsymbol{\zeta _{p}},\boldsymbol{\zeta _{q}})}{F_{\mathbf{Z}%
}^{2}(\mathrm{VaR}_{\boldsymbol{p}}(\mathbf{Z}),\mathrm{VaR}_{\boldsymbol{q}%
}(\mathbf{Z}))}\right] ,
\end{equation*}%
$\mathrm{e}^{\boldsymbol{\mu }}=(\mathrm{e}^{\mu _{1}},\mathrm{e}^{\mu
_{2}},\cdots ,\mathrm{e}^{\mu _{n}})^{\mathrm{T}}$, $\boldsymbol{\psi }%
=\left( \psi_{n} \left( -\frac{1}{2}\sigma _{11}\right) ,\psi_{n} \left( -%
\frac{1}{2}\sigma _{22}\right) ,\cdots ,\psi_{n} \left( -\frac{1}{2}\sigma
_{nn}\right) \right) ^{\mathrm{T}},$ $\boldsymbol{\delta _{\boldsymbol{p,q}}}%
=(\delta _{1},\delta _{2},\cdots ,\delta _{n})^{\mathrm{T}}$, $\delta _{k}=%
\frac{F_{\mathbf{P}_{k}\dagger }(\boldsymbol{\zeta _{p}},\boldsymbol{\zeta
_{q}})}{F_{\mathbf{Z}}(\mathrm{VaR}_{\boldsymbol{p}}(\mathbf{Z}),\mathrm{VaR}%
_{\boldsymbol{q}}(\mathbf{Z}))}$, $k\in \{1,2,\cdots ,n\}$, $\boldsymbol{%
\zeta _{v}}=\mathbf{\Sigma }^{-1/2}[\mathrm{ln}(\mathrm{VaR}_{\boldsymbol{v}%
}(\mathbf{Z}))-\boldsymbol{\mu }]$, $\boldsymbol{v}\in \{\boldsymbol{p},%
\boldsymbol{q}\}$, and symbol $\circ $ is the Hadamard product.
\end{proposition}

Note that for some special distributions, the characteristic generator $%
\psi_{n} $ can take different values, such as normal distribution $\psi_{n}
(u)=\mathrm{e}^{-u}$, Laplace distribution $\psi_{n} (u)=\frac{1}{1+u}$, and
generalized stable laws distributions $\psi_{n} (u)=\mathrm{e}^{-ru^{\alpha
/2}},~\alpha ,r>0 $, etc.

\section{Numerical Illustrations\label{sec:7.1}}

In this section, we present numerical examples to illustrate the proposed
multivariate range risk measures. We first compare the MRVaR and the MRCov
(MRCorr) of the normal, student-$t$ and Laplace distributions for different
dependence structures by using the explicit formulas. Then we consider an
application of optimal portfolio selection. More specifically, we propose a
range-based mean-variance framework of the optimal portfolio selection.

\subsection{A comparison of multivariate range-based risk measures in the
elliptical family}

We consider three triplet random vectors of risks, namely, $\mathbf{U}$, $%
\mathbf{V}$ and $\mathbf{W}$, which follow the normal, Student-$t$ and
Laplace distributions in the elliptical family respectively; i.e., $\mathbf{U%
}=(U_{1},U_{2},U_{3})^{T}\sim N_{3}(\boldsymbol{\mu },\mathbf{\Sigma }),$ $%
\mathbf{V}=(V_{1},V_{2},V_{3})^{T}\sim St_{3}(\boldsymbol{\mu },\mathbf{%
\Sigma },m=4),$ and $\mathbf{W}=(W_{1},W_{2},W_{3})^{T}\sim La_{3}(%
\boldsymbol{\mu },\mathbf{\Sigma })$, we consider vector of means $%
\boldsymbol{\mu }$ and scale matrix $\mathbf{\Sigma }$, respectively,
\begin{equation*}
\boldsymbol{\mu }=\left(
\begin{array}{c}
1.4 \\
1.1 \\
3.4%
\end{array}%
\right) ,~\mathbf{\Sigma }=\left(
\begin{array}{ccc}
1.33 & -0.067 & 0.83 \\
-0.067 & 0.25 & -0.50 \\
0.83 & -0.50 & 5.76%
\end{array}%
\right) .
\end{equation*}%
Let $\boldsymbol{p_{1}}=(0,~0,~0)^{T}$,~$\boldsymbol{q_{1}}%
=(0.10,~0.10,~0.10)^{T},$ $\boldsymbol{p_{2}}=(0.30,~0.30,~0.30)^{T}$, ~$%
\boldsymbol{q_{2}}=(0.70,~0.70,~0.70)^{T},$ \newline
$\boldsymbol{p_{3}}=(0.30,~0.30,~0.30)^{T}$, $\boldsymbol{q_{3}}%
=(0.80,~0.80,~0.80)^{T},$ $\boldsymbol{p_{4}}=(0.95,~0.95,~0.95)^{T}$ and $%
\boldsymbol{q_{4}}=(1,~1,~1)^{T}$. Firstly, we compute VaRs of those
distributions; see Table 1.

\begin{table}[tbh]
\centering
{ 
\centering
 {~~~~~~~~~~~~~~~Table 1: The VaRs of $U_{1}$, $U_{2}$, $U_{3}$, $V_{1}$, $V_{2}$,
$V_{3}$, $W_{1}$, $W_{2}$, $W_{3}$ for different $p$.}\newline
\scriptsize
\begin{tabular}{cccccccccc}
\hline\hline
& $U_{1}$ & $U_{2}$ & $U_{3}$ & $V_{1}$ & $V_{2}$ & $V_{3}$ & $W_{1}$ & $%
W_{2}$ & $W_{3}$ \\ \hline
0 & $-\infty$ & $-\infty$ & $-\infty$ & $-\infty$ & $-\infty$ & $-\infty$ & $%
-\infty$ & $-\infty$ & $-\infty$ \\
0.10 & -0.077958 & 0.459224 & 0.324277 & -0.368180 & 0.333397 & -0.279695 &
-1.364674 & -0.098638 & -2.353460 \\
0.30 & 0.795232 & 0.837800 & 2.141439 & 0.744202 & 0.815675 & 2.035242 &
0.393145 & 0.663478 & 1.304670 \\
0.70 & 2.004768 & 1.362200 & 4.658561 & 2.055799 & 1.384325 & 4.764757 &
2.406857 & 1.536528 & 5.495330 \\
0.80 & 2.370604 & 1.520811 & 5.419894 & 2.485174 & 1.570482 & 5.658313 &
3.090910 & 1.833102 & 6.918880 \\
0.95 & 3.296937 & 1.922426 & 7.347650 & 3.858566 & 2.165923 & 8.516429 &
5.173240 & 2.735906 & 11.252350 \\
1.00 & $\infty$ & $\infty$ & $\infty$ & $\infty$ & $\infty$ & $\infty$ & $%
\infty$ & $\infty$ & $\infty$ \\ \hline\hline
\end{tabular}%
}
\end{table}

Next, we calculate the MRVaRs and the MRCovs of $\mathbf{U}$, $\mathbf{V}$
and $\mathbf{W}$ for different $(\boldsymbol{p,q})$ as follows:

\begin{table}[tbh]
\centering {~~~~~~~~~~~~~~~~~~~~~~~~~~~~Table 2: The MRVaRs of $\mathbf{U}$ for different $(\boldsymbol{%
p,q})$.}\newline
\begin{tabular}{cccc}
\hline\hline
& $U_{1}$ & $U_{2}$ & $U_{3}$ \\ \hline
$(\boldsymbol{p_{1},q_{1}})$ & -0.702948 & 0.296865 & -0.553282 \\
$(\boldsymbol{p_{2},q_{2}})$ & 1.4 & 1.1 & 3.4 \\
$(\boldsymbol{p_{3},q_{3}})$ & 1.561117 & 1.161038 & 3.704342 \\
$(\boldsymbol{p_{4},q_{4}})$ & 3.858560 & 2.061853 & 8.104730 \\ \hline\hline
\end{tabular}%
\end{table}
\begin{equation*}
\mathrm{MRCov}_{(\boldsymbol{p_{1}},\boldsymbol{q_{1}})}(\mathbf{U})=\left(
\begin{array}{ccc}
0.2725561 & -6.3955180\times 10^{-5} & 0.0210362 \\
-6.3955180\times 10^{-5} & 0.02197306 & -0.0053686 \\
0.0210362 & -0.0053686 & 0.6122263%
\end{array}%
\right) ,
\end{equation*}%
\begin{equation*}
\mathrm{MRCov}_{(\boldsymbol{p_{2}},\boldsymbol{q_{2}})}(\mathbf{U})=\left(
\begin{array}{ccc}
0.1171168 & -1.4123580\times 10^{-5} & 0.0071105 \\
-1.4123580\times 10^{-5} & 0.02193414 & -0.0046652 \\
0.0071105 & -0.0046652 & 0.5036470%
\end{array}%
\right) ,
\end{equation*}%
\begin{equation*}
\mathrm{MRCov}_{(\boldsymbol{p_{3}},\boldsymbol{q_{3}})}(\mathbf{U})=\left(
\begin{array}{ccc}
0.1929618 & -0.0001503 & 0.0191361 \\
-0.0001503 & 0.0359268 & -0.0124602 \\
0.0191361 & -0.0124602 & 0.8243330%
\end{array}%
\right) ,
\end{equation*}%
\begin{equation*}
\mathrm{MRCov}_{(\boldsymbol{p_{4}},\boldsymbol{q_{4}})}(\mathbf{U})=\left(
\begin{array}{ccc}
0.2310622 & -1.5509710\times 10^{-5} & 0.0138496 \\
-1.5509710\times 10^{-5} & 0.0168436 & -0.0032329 \\
0.0138496 & -0.0032329 & 0.4767768%
\end{array}%
\right) .
\end{equation*}

\begin{table}[tbh]
\centering ~~~~~~~~~~~~~~~~~~~~~~~~~~~Table 3: The MRVaRs of $\mathbf{V}$ for different $(\boldsymbol{%
p,q})$.\newline
\begin{tabular}{cccc}
\hline\hline
& $V_{1}$ & $V_{2}$ & $V_{3}$ \\ \hline
$(\boldsymbol{p_{1},q_{1}})$ & -2.703324 & -0.355314 & -3.940607 \\
$(\boldsymbol{p_{2},q_{2}})$ & 1.400001 & 1.1 & 3.4 \\
$(\boldsymbol{p_{3},q_{3}})$ & 1.568532 & 1.159554 & 3.705519 \\
$(\boldsymbol{p_{4},q_{4}})$ & 6.867393 & 3.048464 & 13.211250 \\
\hline\hline
\end{tabular}%
\end{table}
\begin{equation*}
\mathrm{MRCov}_{(\boldsymbol{p_{1}},\boldsymbol{q_{1}})}(\mathbf{V})=\left(
\begin{array}{ccc}
7.5413450 & 0.8690497 & 5.3625440 \\
0.8690497 & 0.8078437 & 1.2198950 \\
5.3625440 & 1.2198950 & 21.5812600%
\end{array}%
\right) ,
\end{equation*}%
\begin{equation*}
\mathrm{MRCov}_{(\boldsymbol{p_{2}},\boldsymbol{q_{2}})}(\mathbf{V})=\left(
\begin{array}{ccc}
0.1329993 & -0.0001762 & 0.0140416 \\
-0.0001762 & 0.0248552 & -0.0091530 \\
0.0140416 & -0.0091530 & 0.5695088%
\end{array}%
\right) ,
\end{equation*}%
\begin{equation*}
\mathrm{MRCov}_{(\boldsymbol{p_{3}},\boldsymbol{q_{3}})}(\mathbf{V})=\left(
\begin{array}{ccc}
0.2224948 & -0.0004433 & 0.0365105 \\
-0.0004433 & 0.0411555 & -0.0216376 \\
0.0365105 & -0.0216376 & 0.9434341%
\end{array}%
\right) ,
\end{equation*}%
\begin{equation*}
\mathrm{MRCov}_{(\boldsymbol{p_{4}},\boldsymbol{q_{4}})}(\mathbf{V})=\left(
\begin{array}{ccc}
12.824010 & 1.494721 & 9.180549 \\
1.494721 & 1.365656 & 2.109936 \\
9.180549 & 2.109936 & 36.532980%
\end{array}%
\right) .
\end{equation*}

\begin{table}[tbh]
\centering ~~~~~~~~~~~~~~~~~~~~~~~~~~~~Table 4: The MRVaRs of $\mathbf{W}$ for different $(\boldsymbol{%
p,q})$.\newline
\begin{tabular}{cccc}
\hline\hline
& $W_{1}$ & $W_{2}$ & $W_{3}$ \\ \hline
$(\boldsymbol{p_{1},q_{1}})$ & -3.448572 & -0.679482 & -5.466375 \\
$(\boldsymbol{p_{2},q_{2}})$ & 1.4 & 1.100002 & 3.400001 \\
$(\boldsymbol{p_{3},q_{3}})$ & 1.641988 & 1.184636 & 3.837478 \\
$(\boldsymbol{p_{4},q_{4}})$ & 7.321639 & 3.316232 & 14.374700 \\
\hline\hline
\end{tabular}%
\end{table}
\begin{equation*}
\mathrm{MRCov}_{(\boldsymbol{p_{1}},\boldsymbol{q_{1}})}(\mathbf{W})=\left(
\begin{array}{ccc}
3.2616840 & 0.0868712 & 0.8060599 \\
0.0868712 & 0.3085260 & 0.0466691 \\
0.8060599 & 0.0466691 & 8.3800210%
\end{array}%
\right) ,
\end{equation*}%
\begin{equation*}
\mathrm{MRCov}_{(\boldsymbol{p_{2}},\boldsymbol{q_{2}})}(\mathbf{W})=\left(
\begin{array}{ccc}
0.2955763 & -0.0012197 & 0.0434348 \\
-0.0012197 & 0.0551933 & -0.0279061 \\
0.0434348 & -0.0279061 & 1.2635410%
\end{array}%
\right) ,
\end{equation*}%
\begin{equation*}
\mathrm{MRCov}_{(\boldsymbol{p_{3}},\boldsymbol{q_{3}})}(\mathbf{W})=\left(
\begin{array}{ccc}
0.5058757 & -0.0001401 & 0.1032881 \\
-0.0001401 & 0.0933769 & -0.0514781 \\
0.1032881 & -0.0514781 & 2.1412330%
\end{array}%
\right) ,
\end{equation*}%
\begin{equation*}
\mathrm{MRCov}_{(\boldsymbol{p_{4}},\boldsymbol{q_{4}})}(\mathbf{W})=\left(
\begin{array}{ccc}
3.5062380 & 0.0777070 & 0.7360139 \\
0.0777070 & 0.3113243 & 0.0393784 \\
0.7360139 & 0.0393784 & 8.5506770%
\end{array}%
\right) .
\end{equation*}

Tables 2-4 present that the MRVaRs of normal, Student-$t$ and Laplace
distribution for different $(\boldsymbol{p,q})$, respectively. We observe
that the MRVaR is increasing in ($\boldsymbol{p_{k}},\boldsymbol{q_{k}}$).
This reveals the more right tail may has larger MRVaR. In particular, at ($%
\boldsymbol{p_{2}},\boldsymbol{q_{2}}$), the MRVaRs of several distributions
are (almost) equal to their expectation $\boldsymbol{\mu }$. This is
consistent with ``the integral of symmetric distribution on symmetric
interval equals to expectation".

As for the $\mathrm{MRCov}$, we find that for Student-$t$ and Laplace
distributions, the diagonal of $\mathrm{MRCov}$ matrix at $(\boldsymbol{%
p_{1},q_{1}})$ and $(\boldsymbol{p_{4},q_{4}})$ is larger than the one at $(%
\boldsymbol{p_{2},q_{2}})$ and $(\boldsymbol{p_{3},q_{3}})$. This means that
for Student-t and Laplace distributions, the tails have larger RVs.
Furthermore, at $(\boldsymbol{p_{1},q_{1}})$ and $(\boldsymbol{p_{4},q_{4}})$%
, the diagonal of $\mathrm{MRCov}$ matrix under Student-$t$ is largest among
all models. By contrast, at $(\boldsymbol{p_{2},q_{2}})$ and $(\boldsymbol{%
p_{3},q_{3}})$, the diagonal of $\mathrm{MRCov}$ matrix under Laplace model
is largest. This indicates that Student-$t$ model may be more heavy-tailed
than Laplace model.

Using Equation (\ref{(v2)}), we can also compute MRCorr (see Appendix). We
can find that for Student-$t$ ($\mathbf{V}$) and Laplace ($\mathbf{W}$)
distributions, the upper triangular entries of $\mathrm{MRCorr}$ matrix at $(%
\boldsymbol{p_{1},q_{1}})$ and $(\boldsymbol{p_{4},q_{4}})$ are larger than
the ones at $(\boldsymbol{p_{2},q_{2}})$ and $(\boldsymbol{p_{3},q_{3}})$
but this does not hold for the $\mathrm{MRCorr}$ matrix of normal ($\mathbf{U%
}$) model. As such, Student-t and Laplace distributions have stronger
pairwise dependence for random vectors. Moreover, $\mathrm{MRCov(\mathbf{V})}%
_{23}$ and $\mathrm{MRCov(\mathbf{W})}_{23}$ are negative at $(\boldsymbol{%
p_{2},q_{2}})$ and $(\boldsymbol{p_{3},q_{3}})$, but turn positive at $(%
\boldsymbol{p_{1},q_{1}})$ and $(\boldsymbol{p_{4},q_{4}})$. Hence, the
chosen range has direct impact on the pairwise dependence of the random
vectors.

Fig.s 1-4 display the MRVaR of the three models at the four chosen ranges
respectively. We may observe that the MRVaR of normal distribution is
largest among the three at $(\boldsymbol{p_{1},q_{1}})$ but turn to the
least at $(\boldsymbol{p_{4},q_{4}})$. Moreover, at $(\boldsymbol{p_{2},q_{2}%
})$ and $(\boldsymbol{p_{3},q_{3}})$, the MRVaR of the three models
distinguish little from each other. In another word, the underlying
distribution may not change the MVaR significantly if we dwell in non-tail
regions whilst the impact extreme event (at tails) are fairly reflected in
MVaR. Therefore, as aforementioned, our range-based risk measures are indeed
very appropriate for practical applications as they suit both investors' and
regulators' appetites.

\subsection{The range-based mean-variance framework of portfolio selection
and efficient frontiers\label{sec:7.2}}

The celebrated Markowitz's mean-variance portfolio theory is the foundation
of the modern portfolio theory. In practice, it is often a certain region of
the underlying risky assets are interested and concerns, especially from the
investors and the shareholders' standpoints. With the explicit formulas at
hand, we can develop a range-based mean-variance framework of the optimal
portfolio selection by directly replacing the means and variances with the
MRVaR and the MRCov. More specifically, for a random vector of $n$%
-dimensional risks $\mathbf{X}=(X_{1},X_{2},\cdots ,X_{n})^{T}$, we consider
an optimization as follows.
\begin{equation}
\min_{\boldsymbol{\omega }}\boldsymbol{\omega }^{T}\mathrm{MRCov}_{(%
\boldsymbol{p,q})}(\mathbf{X})\boldsymbol{\omega },  \label{(v16)}
\end{equation}%
subject to the portfolio constrain $\boldsymbol{1}^{T}\boldsymbol{\omega }=1$
and a prespecified target level $\mathrm{MRVaR}_{(\boldsymbol{p,q})}(\mathbf{%
X})^{T}\boldsymbol{\omega }=\mu _{0}$, where $\boldsymbol{\omega }=(\omega
_{1},\cdots ,\omega _{n})^{T}$ and $\boldsymbol{1}=(1,\cdots ,1)^{T}$. As an
analogue to the mean-variance portfolio theory, the following proposition
gives the efficient portfolio that solves the optimization.

\begin{proposition}
\label{pro.4} Assume $\mathrm{MRCov}_{(\boldsymbol{p,q})}(\mathbf{X})$ is
invertible, the unique optimal solution to (\ref{(v16)}) is
\begin{align*}
\boldsymbol{\omega }^{\ast }=& \frac{1}{d}\bigg\{\left[ a\mathrm{MRCov}_{(%
\boldsymbol{p,q})}(\mathbf{X})^{-1}\boldsymbol{1}-b\mathrm{MRCov}_{(%
\boldsymbol{p,q})}(\mathbf{X})^{-1}\mathrm{MRVaR}_{(\boldsymbol{p,q})}(%
\mathbf{X})\right]  \\
& +\mu _{0}\left[ c\mathrm{MRCov}_{(\boldsymbol{p,q})}(\mathbf{X})^{-1}%
\mathrm{MRVaR}_{(\boldsymbol{p,q})}(\mathbf{X})-b\mathrm{MRCov}_{(%
\boldsymbol{p,q})}(\mathbf{X})^{-1}\boldsymbol{1}\right] \bigg\},
\end{align*}%
where\newline
\begin{equation*}
a=\mathrm{MRVaR}_{(\boldsymbol{p,q})}(\mathbf{X})^{T}\mathrm{MRCov}_{(%
\boldsymbol{p,q})}(\mathbf{X})^{-1}\mathrm{MRVaR}_{(\boldsymbol{p,q})}(%
\mathbf{X}),~b=\boldsymbol{1}^{T}\mathrm{MRCov}_{(\boldsymbol{p,q})}(\mathbf{%
X})^{-1}\mathrm{MRVaR}_{(\boldsymbol{p,q})}(\mathbf{X}),
\end{equation*}%
\begin{equation*}
~c=\boldsymbol{1}^{T}\mathrm{MRCov}_{(\boldsymbol{p,q})}(\mathbf{X})^{-1}%
\boldsymbol{1}~\mathrm{and}~d=ac-b^{2}.
\end{equation*}
\end{proposition}

\noindent

We further consider an example based on stock daily log-return data. We take
five stocks from the Nasdaq stock market (Apple Inc. (AAPL), Cisco Sys. Inc.
(CSCO), eBay Inc. (EBAY), Intel Corporation (INTC), Sirius XM Holdings Inc.
(SIRI)) from the year 2020 to 2021. We use student-$t$ distribution to fit
data, and estimate the parameters via maximum likelihood estimation. We
denote the random returns of the stocks as $\mathbf{X}\sim St_{5}(%
\boldsymbol{\mu },\mathbf{\Sigma },6.2623761)$, where
\begin{equation*}
\boldsymbol{\mu }=10^{-3}\left(
\begin{array}{c}
1.1216 \\
0.5261 \\
0.1270 \\
-0.6643 \\
-0.2408%
\end{array}%
\right) ,\mathbf{\Sigma }=10^{-2}\left(
\begin{array}{ccccc}
1.12591 & 0.28224 & 0.25971 & 0.34555 & 0.13225 \\
0.28224 & 0.95803 & 0.13244 & 0.20053 & 0.14266 \\
0.25971 & 0.13244 & 1.58482 & 0.21797 & 0.11606 \\
0.34555 & 0.20053 & 0.21797 & 1.14941 & 0.18298 \\
0.13225 & 0.14266 & 0.11606 & 0.18298 & 1.16874%
\end{array}%
\right) .
\end{equation*}%
We consider three chosen ranges\footnote{%
Note that we can replace the VaR at $(\boldsymbol{p,q})$ in the optimal
portfolio weights by simple thresholds, i.e. $k_{i}\leq X_{i}\leq l_{i},$%
where $k_{i}$ and $l_{i}$ are chosen constants, $i=1,2,...,n.$}, $%
\boldsymbol{p_{1}}=(0.75,~0.75,~0.75,~0.75,~0.75)^{T}$,~$\boldsymbol{q_{1}}%
=(0.95,~0.95,~0.95,~0.95,~0.95)^{T},$ $\boldsymbol{p_{2}}%
=(0.80,~0.80,~0.80,~0.80,~0.80)^{T}$, $\boldsymbol{q_{2}}=(1,~1,~1,~1,~1)^{T}
$, $\boldsymbol{p_{3}}=(0.95,~0.95,~0.95,~0.95,~0.95)^{T}$ and $\boldsymbol{%
q_{3}}=(1,~1,~1,~1,~1)^{T},$ then the VaRs of $X_{1}$, $X_{2}$, $X_{3}$, $%
X_{4}$ and $X_{5}$ are shown in Table 5 and the MRVaRs of $\mathbf{X}$ are
presented in Table 6.
\begin{table}[tbh]
\centering ~~~~~~~~~~~~~~~~~~~~~~Table 5: The VaRs of $X_{1}$, $X_{2}$, $X_{3}$, $X_{4}$ and $%
X_{5} $ for different $p$.\newline
\begin{tabular}{cccccc}
\hline\hline
& $X_{1}$ & $X_{2}$ & $X_{3}$ & $X_{4}$ & $X_{5}$ \\ \hline
0.40 & -0.026928 & -0.025348 & -0.033151 & -0.029005 & -0.028818 \\
0.60 & 0.029171 & 0.026400 & 0.033405 & 0.027676 & 0.028337 \\
0.70 & 0.059708 & 0.054568 & 0.069634 & 0.058529 & 0.059449 \\
0.75 & 0.077060 & 0.070575 & 0.090222 & 0.076063 & 0.077129 \\
0.80 & 0.096922 & 0.088896 & 0.113786 & 0.096131 & 0.097365 \\
0.90 & 0.153112 & 0.140727 & 0.180451 & 0.152904 & 0.154613 \\
0.95 & 0.205775 & 0.189307 & 0.242932 & 0.206114 & 0.208267 \\
1.00 & $\infty$ & $\infty$ & $\infty$ & $\infty$ & $\infty$ \\ \hline\hline
\end{tabular}%
\end{table}
\begin{table}[tbh]
\centering
{
\centering
~~~~~~~~~~~~~~~~~~~~Table 6: The MRVaRs of $\mathbf{X}$ for $(\boldsymbol{p_{1}},\boldsymbol{%
q_{1}}),~(\boldsymbol{p_{2}},\boldsymbol{q_{2}})~ \mathrm{and} ~(\boldsymbol{%
p_{3}},\boldsymbol{q_{3}})$.\newline
\begin{tabular}{|c||c|c|c|c|c|}
\hline
\backslashbox
{$(\boldsymbol{p_{i}},\boldsymbol{q_{i}})$}{$\mathrm{MRVaR}_{(%
\boldsymbol{p_{i}},\boldsymbol{q_{i}})}(\mathbf{X})$}{$X_{i}$} & $X_{1}$ & $%
X_{2}$ & $X_{3}$ & $X_{4}$ & $X_{5}$ \\ \hline
$(\boldsymbol{p_{1}},\boldsymbol{q_{1}})$ & 0.1332060 & 0.1206013 & 0.1530140
& 0.1320112 & 0.1307101 \\ \hline
$(\boldsymbol{p_{2}},\boldsymbol{q_{2}})$ & 0.2371618 & 0.2091944 & 0.2603811
& 0.2338855 & 0.2213374 \\ \hline
$(\boldsymbol{p_{3}},\boldsymbol{q_{3}})$ & 0.4030022 & 0.3569983 & 0.4461190
& 0.3991542 & 0.3798648 \\ \hline
\end{tabular}
}
\end{table}
Moreover, the MRCovs of $\mathbf{X}$, by Corollary \ref{co.7}, are
\begin{equation*}
\mathrm{MRCov}_{(\boldsymbol{p_{1}},\boldsymbol{q_{1}})}(\mathbf{X}%
)=10^{-3}\left(
\begin{array}{ccccc}
1.27195500 & 0.04475777 & 0.04455139 & 0.05303480 & 0.02447907 \\
0.04475777 & 1.07132200 & 0.02831541 & 0.03035051 & 0.03088864 \\
0.04455139 & 0.02831541 & 1.75647200 & 0.03894598 & 0.03414841 \\
0.05303480 & 0.03035051 & 0.03894598 & 1.29421400 & 0.03557834 \\
0.02447907 & 0.03088864 & 0.03414841 & 0.03557834 & 1.29211600%
\end{array}%
\right) ,
\end{equation*}%
\begin{equation*}
\mathrm{MRCov}_{(\boldsymbol{p_{2}},\boldsymbol{q_{2}})}(\mathbf{X}%
)=10^{-2}\left(
\begin{array}{ccccc}
1.6051150 & 0.4483114 & 0.5045233 & 0.5268113 & 0.3850687 \\
0.4483114 & 1.2609090 & 0.3945301 & 0.3959148 & 0.3503502 \\
0.5045233 & 0.3945301 & 1.9562300 & 0.4777596 & 0.4089094 \\
0.5268113 & 0.3959148 & 0.4777596 & 1.5922910 & 0.4076230 \\
0.3850687 & 0.3503502 & 0.4089094 & 0.4076230 & 1.4156880%
\end{array}%
\right)
\end{equation*}%
and
\begin{equation*}
\mathrm{MRCov}_{(\boldsymbol{p_{3}},\boldsymbol{q_{3}})}(\mathbf{X}%
)=10^{-2}\left(
\begin{array}{ccccc}
3.4845110 & 1.0072270 & 1.1441760 & 1.1812240 & 0.8861056 \\
1.0072270 & 2.7075940 & 0.9071066 & 0.9011320 & 0.7995451 \\
1.1441760 & 0.9071066 & 4.1678580 & 1.0889430 & 0.9398218 \\
1.1812240 & 0.9011320 & 1.0889430 & 3.4419590 & 0.9280654 \\
0.8861056 & 0.7995451 & 0.9398218 & 0.9280654 & 3.0099370%
\end{array}%
\right) .
\end{equation*}%
We then calculate the efficient frontier of the optimal portfolios under the
three chosen ranges (see Fig.5). We can observe that including the tails can
result in very different optimal portfolios, comparing with excluding tails.
Although it is generally fair that optimal portfolio on the upper tails
region has higher expected return, the optimal portfolio excluding the tails
may outperform the ones that only concern the tail ranges, which is again
consistent with our former discussions that the shareholders and investors
may be more interested in the non-tail ranges.

\section{Concluding remarks\label{sec:8}}

In this paper, we propose multivariate range-based risk measures, namely the
MRVaR and MRCov (MRCor). Both extend the univariate corresponding risk
measures in the literature and have much more flexibility in choosing the
concerned ranges, which is of importance and practical motivations. In the
context of multivariate elliptical distribution, we derive explicit formulas
of MRVaR and MRCov and use them for analytical analysis. In particular, we
propose a range-based mean-variance framework that are more flexible in
optimal portfolio selection. By deriving the range-based efficient
frontiers, we also show that such optimal portfolio selection may serve for
more general preferences in practice.

\section*{Acknowledgments}

\noindent The research was supported by the National Natural Science
Foundation of China (No. 12071251)

\section*{Conflicts of Interest}

\noindent The authors declare that they have no conflicts of interest.

\section*{Appendix}

\noindent$\mathbf{Proof ~of~ Proposition~ 1.}$ (i) Using definition of the
MRVaR and the positive homogeneity property of $\mathrm{VaR}$, we directly
have
\begin{align*}
\mathrm{MRVaR}_{(\boldsymbol{p},\boldsymbol{q})}(c\mathbf{X})&=\mathrm{E}%
\left[c\mathbf{X}|c\mathrm{VaR}_{\boldsymbol{p}}(\mathbf{X})\leq c\mathbf{X}%
\leq c\mathrm{VaR}_{\boldsymbol{q}}(\mathbf{X})\right] \\
&=c\mathrm{E}\left[\mathbf{X}|\mathrm{VaR}_{\boldsymbol{p}}(\mathbf{X})\leq%
\mathbf{X}\leq\mathrm{VaR}_{\boldsymbol{q}}(\mathbf{X})\right].
\end{align*}
(ii) Using definition of the MRVaR, translation invariance of VaR and the
linear property of expectation, we directly have
\begin{align*}
\mathrm{MRVaR}_{(\boldsymbol{p},\boldsymbol{q})}(\mathbf{X}+\boldsymbol{%
\gamma})&=\mathrm{E}\left[\mathbf{X}+\boldsymbol{\gamma}|\mathrm{VaR}_{%
\boldsymbol{p}}(\mathbf{X})+\boldsymbol{\gamma}\leq\mathbf{X}+\boldsymbol{%
\gamma}\leq\mathrm{VaR}_{\boldsymbol{q}}(\mathbf{X})+\boldsymbol{\gamma}%
\right] \\
&=\boldsymbol{\gamma}+\mathrm{E}\left[\mathbf{X}|\mathrm{VaR}_{\boldsymbol{p}%
}(\mathbf{X})\leq\mathbf{X}\leq\mathrm{VaR}_{\boldsymbol{q}}(\mathbf{X})%
\right].
\end{align*}
(iii) If $X_{1},~X_{2},\cdots,X_{n}$ are independent, then
\begin{align*}
\mathrm{MRVaR}_{(\boldsymbol{p},\boldsymbol{q})}(\mathbf{X})&=\left(%
\begin{array}{c}
\mathrm{E}\left(X_{1}|\mathrm{VaR}_{\boldsymbol{p}}(\mathbf{X})\leq\mathbf{X}%
\leq \mathrm{VaR}_{\boldsymbol{q}}(\mathbf{X})\right) \\
\mathrm{E}\left(X_{2}|\mathrm{VaR}_{\boldsymbol{p}}(\mathbf{X})\leq\mathbf{X}%
\leq \mathrm{VaR}_{\boldsymbol{q}}(\mathbf{X})\right) \\
\cdots \\
\mathrm{E}\left(X_{n}|\mathrm{VaR}_{\boldsymbol{p}}(\mathbf{X})\leq\mathbf{X}%
\leq \mathrm{VaR}_{\boldsymbol{q}}(\mathbf{X})\right)%
\end{array}
\right) \\
&=\left(%
\begin{array}{c}
\mathrm{E}\left(X_{1}|\mathrm{VaR}_{p_{1}}(X_{1})\leq X_{1}\leq \mathrm{VaR}%
_{q_{1}}(X_{1})\right) \\
\mathrm{E}\left(X_{2}|\mathrm{VaR}_{p_{2}}(X_{2})\leq X_{2}\leq \mathrm{VaR}%
_{q_{2}}(X_{2})\right) \\
\cdots \\
\mathrm{E}\left(X_{n}|\mathrm{VaR}_{p_{n}}(X_{n})\leq X_{n}\leq \mathrm{VaR}%
_{q_{n}}(X_{n})\right)%
\end{array}
\right) =\left(%
\begin{array}{c}
\mathrm{RVaR}_{(p_{1},q_{1})}(X_{1}) \\
\mathrm{RVaR}_{(p_{2},q_{2})}(X_{2}) \\
\cdots \\
\mathrm{RVaR}_{(p_{n},q_{n})}(X_{n})%
\end{array}
\right).
\end{align*}
(iv) Since $\mathbf{Y}\geq \mathbf{X}$, we have $\mathbf{W}=\mathbf{Y}-%
\mathbf{X}\geq\boldsymbol{0}$ is a non-negative random vector, thus for any $%
\boldsymbol{p}<\boldsymbol{q}$, $\boldsymbol{p,q}\in(0,1)^{n}$, $\mathrm{VaR}%
_{\boldsymbol{q}}(\mathbf{W})>\mathrm{VaR}_{\boldsymbol{p}}(\mathbf{W})\geq%
\boldsymbol{0}$. Therefore, we can write the $\mathrm{MRE}_{(\boldsymbol{p},%
\boldsymbol{\ q})}(\mathbf{W})$,
\begin{align*}
\mathrm{MRVaR}_{(\boldsymbol{p},\boldsymbol{q})}(\mathbf{W})&=\mathrm{E}%
\left[\mathbf{W}|\mathrm{VaR}_{\boldsymbol{p}}(\mathbf{W})\leq\mathbf{W}\leq%
\mathrm{VaR}_{\boldsymbol{q}}(\mathbf{W})\right],
\end{align*}
and since $\mathrm{VaR}_{\boldsymbol{q}}(\mathbf{W})>\mathrm{VaR}_{%
\boldsymbol{p}}(\mathbf{W})\geq\boldsymbol{0}$ it is clear that $\mathrm{%
MRVaR}_{(\boldsymbol{p},\boldsymbol{q})}(\mathbf{W})\geq\boldsymbol{0}$,
completing the proof.\newline
\noindent$\mathbf{Proof ~of~ Proposition~ 2.}$ (i) Using definition of the
MRCov and $\mathrm{MRVaR}_{(\boldsymbol{p},\boldsymbol{q})}(\boldsymbol{%
\gamma})=\boldsymbol{\gamma}$, we readily obtain (i). \newline
(ii) Using definition of the MRCov and the positive homogeneity property of $%
\mathrm{VaR}$, we have
\begin{align*}
&\mathrm{MRCov}_{(\boldsymbol{p},\boldsymbol{q})}(c\mathbf{X})=\mathrm{E}%
\left[(c\mathbf{X}-\mathrm{MRVaR}_{(\boldsymbol{p},\boldsymbol{q})}(c\mathbf{%
X}))(c\mathbf{X}-\mathrm{MRVaR}_{(\boldsymbol{p},\boldsymbol{q})}(c\mathbf{X}%
))^{\mathrm{T}}|c\mathrm{VaR}_{\boldsymbol{p}}(\mathbf{X})\leq c\mathbf{X}%
\leq c\mathrm{VaR}_{\boldsymbol{q}}(\mathbf{X})\right],
\end{align*}
and from the positive homogeneity property of $\mathrm{MRVaR}$, we obtain
\begin{align*}
\mathrm{MRCov}_{(\boldsymbol{p},\boldsymbol{q})}(c\mathbf{X})&=\mathrm{E}%
\left[(c\mathbf{X}-c\mathrm{MRVaR}_{(\boldsymbol{p},\boldsymbol{q})}(\mathbf{%
X}))(c\mathbf{X}-c\mathrm{MRVaR}_{(\boldsymbol{p},\boldsymbol{q})}(\mathbf{X}%
))^{\mathrm{T}}|c\mathrm{VaR}_{\boldsymbol{p}}(\mathbf{X})\leq c\mathbf{X}%
\leq c\mathrm{VaR}_{\boldsymbol{q}}(\mathbf{X})\right] \\
&=c^{2}\mathrm{MRCov}_{(\boldsymbol{p},\boldsymbol{q})}(\mathbf{X}).
\end{align*}
(iii) Using definition of the MRCov, the linear property of $\mathrm{VaR}$
and Eq. (ii) of MRVaR, we directly have
\begin{align*}
&\mathrm{MRCov}_{(\boldsymbol{p},\boldsymbol{q})}(\mathbf{X}+\boldsymbol{%
\gamma}) \\
&=\mathrm{E}\left[(\mathbf{X}+\boldsymbol{\gamma}-\mathrm{MRVaR}_{(%
\boldsymbol{p},\boldsymbol{q})}(\mathbf{X}+\boldsymbol{\gamma}))(\mathbf{X}+%
\boldsymbol{\gamma}-\mathrm{MRVaR}_{(\boldsymbol{p},\boldsymbol{q})}(\mathbf{%
X}+\boldsymbol{\gamma}))^{\mathrm{T}}|\mathrm{VaR}_{\boldsymbol{p}}(\mathbf{X%
})+\boldsymbol{\gamma}\leq\mathbf{X}+\boldsymbol{\gamma}\leq\mathrm{VaR}_{%
\boldsymbol{q}}(\mathbf{X})+\boldsymbol{\gamma}\right] \\
&=\mathrm{MRCov}_{(\boldsymbol{p},\boldsymbol{q})}(\mathbf{X}).
\end{align*}
(iv) If $X_{1},~X_{2},\cdots,X_{n}$ are independent, then
\begin{align*}
\mathrm{MRCov}_{(\boldsymbol{p},\boldsymbol{q})}(\mathbf{X})&=\mathrm{E}%
\left[(\mathbf{X}-\mathrm{MRVaR}_{(\boldsymbol{p},\boldsymbol{q})}(\mathbf{X}%
))(\mathbf{X}-\mathrm{MRVaR}_{(\boldsymbol{p},\boldsymbol{q})}(\mathbf{X}))^{%
\mathrm{T}}|\mathrm{VaR}_{\boldsymbol{p}}(\mathbf{X})\leq\mathbf{X}\leq%
\mathrm{VaR}_{\boldsymbol{q}}(\mathbf{X})\right] \\
&=\mathrm{E}\left[\mathbf{X}\mathbf{X}^{\mathrm{T}}|\mathrm{VaR}_{%
\boldsymbol{p}}(\mathbf{X})\leq\mathbf{X}\leq\mathrm{VaR}_{\boldsymbol{q}}(%
\mathbf{X})\right]-\mathrm{MRVaR}_{(\boldsymbol{p},\boldsymbol{q})}(\mathbf{X%
})\mathrm{MRVaR}_{(\boldsymbol{p},\boldsymbol{q})}(\mathbf{X})^{\mathrm{T}}
\\
&=\mathrm{E}\left[\mathbf{X}\mathbf{X}^{\mathrm{T}}|\mathrm{VaR}_{%
\boldsymbol{p}}(\mathbf{X})\leq\mathbf{X}\leq\mathrm{VaR}_{\boldsymbol{q}}(%
\mathbf{X})\right]-(\mathrm{RVaR}_{(p_{1},q_{1})}(X_{1}),\mathrm{RVaR}%
_{(p_{2},q_{2})}(X_{2}),\cdots,\mathrm{RVaR}_{(p_{n},q_{n})}(X_{n})) \\
&~~~\boldsymbol{\cdot}(\mathrm{RVaR}_{(p_{1},q_{1})}(X_{1}),\mathrm{RVaR}%
_{(p_{2},q_{2})}(X_{2}),\cdots,\mathrm{RVaR}_{(p_{n},q_{n})}(X_{n}))^{%
\mathrm{T}} \\
&=\mathrm{diag}(\mathrm{RV}_{(p_{1},q_{1})}(X_{1}),\mathrm{RV}%
_{(p_{2},q_{2})}(X_{2}),\cdots,\mathrm{RV}_{(p_{n},q_{n})}(X_{n})),
\end{align*}
as required.

\noindent$\mathbf{Proof ~of~ Theorem~ 1.}$ Since the RVaR is a special case
of MRVaR ($n=1$), the proof of Theorem \ref{th.1.1} is similar to that of
Theorem \ref{th.1}. We will omit it here.

\noindent$\mathbf{Proof ~of~ Theorem~ 2.}$ Using definition, we have
\begin{align*}
&\mathrm{MRVaR}_{(\boldsymbol{p,q})}(\mathbf{X})=\frac{c_{n}}{\sqrt{|\mathbf{%
\Sigma}|}F_{\mathbf{X}}(\boldsymbol{a,b})}\int_{\boldsymbol{a}}^{\boldsymbol{%
b}}\boldsymbol{x}g_{n}\left(\frac{1}{2}\boldsymbol{(x-\mu)^{T}\Sigma^{-1}(x-%
\mu)}\right)\mathrm{d}\boldsymbol{x},
\end{align*}
where $\boldsymbol{a}=VaR_{\boldsymbol{p}}(\mathbf{X})$ and $\boldsymbol{b}%
=VaR_{\boldsymbol{q}}(\mathbf{X})$.\newline
Applying translation $\mathbf{Y}=\mathbf{\Sigma}^{-\frac{1}{2}}(\mathbf{X}-%
\boldsymbol{\mu})$, we obtain
\begin{align*}
\mathrm{MRVaR}_{(\boldsymbol{p,q})}(\mathbf{X})&=\frac{c_{n}}{F_{\mathbf{Y}}(%
\boldsymbol{\eta_{p},\eta_{q}})}\int_{\boldsymbol{\eta_{p}}}^{\boldsymbol{%
\eta_{q}}}\left(\mathbf{\Sigma}^{\frac{1}{2}}\boldsymbol{y}+\boldsymbol{\mu}%
\right) g_{n}\left(\frac{1}{2}\boldsymbol{y^{T}y}\right)\mathrm{d}%
\boldsymbol{y} \\
&=\boldsymbol{\mu}+\frac{\mathbf{\Sigma}^{\frac{1}{2}}}{F_{\mathbf{Y}}(%
\boldsymbol{\eta_{p},\eta_{q}})}\boldsymbol{\delta_{p,q}},
\end{align*}
where
\begin{equation*}
\boldsymbol{\delta_{p,q}}=\int_{\boldsymbol{\eta_{p}}}^{\boldsymbol{\eta_{q}}%
}c_{n}\boldsymbol{y} g_{n}\left(\frac{1}{2}\boldsymbol{y^{T}y}\right)\mathrm{%
d}\boldsymbol{y}.
\end{equation*}
Note that
\begin{align*}
\delta_{k,\boldsymbol{p,q}}&=c_{n}\int_{\boldsymbol{\eta_{p}}}^{\boldsymbol{%
\eta_{q}}}y_{k}g_{n}\left(\frac{1}{2}\boldsymbol{y^{T}y}\right)\mathrm{d}%
\boldsymbol{y} \\
&=c_{n}\int_{\boldsymbol{\eta}_{\boldsymbol{p},-k}}^{\boldsymbol{\eta}_{%
\boldsymbol{q},-k}}\int_{\eta_{\boldsymbol{p},k}}^{\eta_{\boldsymbol{q}%
,k}}-\partial_{k}\overline{G}_{n}\left(\frac{1}{2}\boldsymbol{y}_{-k}^{T}%
\boldsymbol{y}_{-k}+\frac{1}{2}y_{k}^{2}\right)\mathrm{d}\boldsymbol{y}_{-k}
\\
&=c_{n}\int_{\boldsymbol{\eta}_{\boldsymbol{p},-k}}^{\boldsymbol{\eta}_{%
\boldsymbol{q},-k}}\left[\overline{G}_{n}\left(\frac{1}{2}\boldsymbol{y}%
_{-k}^{T}\boldsymbol{y}_{-k}+\frac{1}{2}\eta_{\boldsymbol{p},k}^{2}\right)-%
\overline{G}_{n}\left(\frac{1}{2}\boldsymbol{y}_{-k}^{T}\boldsymbol{y}_{-k}+%
\frac{1}{2}\eta_{\boldsymbol{q},k}^{2}\right)\right]\mathrm{d}\boldsymbol{y}%
_{-k} \\
&=c_{n}\left\{\frac{1}{c_{n-1,\boldsymbol{p},k}^{\ast}}F_{\mathbf{Y}_{%
\boldsymbol{p},-k}}(\boldsymbol{\eta}_{\boldsymbol{p},-k},\boldsymbol{\eta}_{%
\boldsymbol{q},-k})-\frac{1}{c_{n-1,\boldsymbol{q},k}^{\ast}}F_{\mathbf{Y}_{%
\boldsymbol{q},-k}}(\boldsymbol{\eta}_{\boldsymbol{p},-k},\boldsymbol{\eta}_{%
\boldsymbol{q},-k})\right\},
\end{align*}
$k=1,~2,\cdots,n.$ Therefore, we obtain the desired result.

\noindent$\mathbf{Proof ~of~ Theorem~ 3.}$ Since the RV is a special case of
MRCov ($n=1$), the proof of Theorem \ref{th.2.2} is similar to that of
Theorem \ref{th.2}. We will omit it here.

\noindent $\mathbf{Proof~of~ Theorem~ 4.}$ Using definition of MRCov, we
have
\begin{align*}
\mathrm{MRCov}_{(\boldsymbol{p},\boldsymbol{q})}(\mathbf{X})&=\mathrm{E}%
\left[(\mathbf{X}-\mathrm{MRVaR}_{(\boldsymbol{p},\boldsymbol{q})}(\mathbf{X}%
))(\mathbf{X}-\mathrm{MRVaR}_{(\boldsymbol{p},\boldsymbol{q})}(\mathbf{X}%
))^{T}|\mathrm{VaR}_{\boldsymbol{p}}(\mathbf{X})\leq\mathbf{X}\leq\mathrm{VaR%
}_{\boldsymbol{q}}(\mathbf{X})\right] \\
&=\mathrm{E}\left[\mathbf{X}\mathbf{X}^{T}|\mathrm{VaR}_{\boldsymbol{p}}(%
\mathbf{X})\leq\mathbf{X}\leq\mathrm{VaR}_{\boldsymbol{q}}(\mathbf{X})\right]%
-\mathrm{MRVaR}_{(\boldsymbol{p},\boldsymbol{q})}(\mathbf{X})\mathrm{MRVaR}%
_{(\boldsymbol{p},\boldsymbol{q})}^{T}(\mathbf{X}).
\end{align*}
Using the transformation $\mathbf{Y}=\mathbf{\Sigma}^{-\frac{1}{2}}(\mathbf{Y%
}-\boldsymbol{\mu})$ and basic algebraic calculations, we obtain
\begin{align*}
&\mathrm{MRCov}_{(\boldsymbol{p},\boldsymbol{q})}(\mathbf{X})=\mathbf{\Sigma}%
^{\frac{1}{2}}\left\{\mathrm{E}[\boldsymbol{Y}\boldsymbol{Y}^{T}|\boldsymbol{%
\eta_{p}}\leq\mathbf{Y}\leq\boldsymbol{\eta_{q}}]-\mathrm{MRVaR}_{(%
\boldsymbol{\eta_{p}},\boldsymbol{\eta_{q}})}(\mathbf{Y})\mathrm{MRVaR}_{(%
\boldsymbol{\eta_{p}},\boldsymbol{\eta_{q}})}^{T}(\mathbf{Y})\right\}\mathbf{%
\Sigma}^{\frac{1}{2}},
\end{align*}
where $\boldsymbol{\eta_{v}}=\mathbf{\Sigma}^{-\frac{1}{2}}(\mathrm{VaR}_{%
\boldsymbol{v}}(X)-\boldsymbol{\mu})$,~$\boldsymbol{v}\in\{\boldsymbol{p,q}%
\} $. \newline
Note that
\begin{align*}
\mathrm{E}[Y_{i}Y_{j}|\boldsymbol{\eta_{p}}\leq\mathbf{Y}\leq\boldsymbol{%
\eta_{q}}]&=\frac{1}{F_{\mathbf{Y}}(\boldsymbol{\eta_{p}},\boldsymbol{%
\eta_{q}})}\int_{\boldsymbol{\eta_{p}}}^{\boldsymbol{\eta_{q}}%
}y_{i}y_{j}c_{n}g_{n}\left(\frac{1}{2}\boldsymbol{y}^{T}\boldsymbol{y}\right)%
\mathrm{d}\boldsymbol{y} \\
&=\frac{c_{n}}{F_{\mathbf{Y}}(\boldsymbol{\eta_{p}},\boldsymbol{\eta_{q}})}%
\int_{\boldsymbol{\eta}_{\boldsymbol{p},-i}}^{\boldsymbol{\eta}_{\boldsymbol{%
q},-i}}y_{j}\int_{\eta_{\boldsymbol{p},i}}^{\eta_{\boldsymbol{q}%
,i}}y_{i}g_{n}\left(\frac{1}{2}\boldsymbol{y}_{-i}^{T}\boldsymbol{y}_{-i}+%
\frac{1}{2}y_{i}^{2}\right)\mathrm{d}y_{i}\mathrm{d}\boldsymbol{y}_{-i},~for
~i\neq j,
\end{align*}
where $\boldsymbol{y}_{-i}=(y_{1},\cdots,y_{i-1},y_{i+1},\cdots,y_{n})^{T}$.%
\newline
Since
\begin{align*}
\int_{\eta_{\boldsymbol{p},i}}^{\eta_{\boldsymbol{q},i}}y_{i}g_{n}\left(%
\frac{1}{2}\boldsymbol{y}_{-i}^{T}\boldsymbol{y}_{-i}+\frac{1}{2}%
y_{i}^{2}\right)\mathrm{d}y_{i}=&-\int_{\eta_{\boldsymbol{p},i}}^{\eta_{%
\boldsymbol{q},i}}\partial_{i}\overline{G}_{n}\left(\frac{1}{2}\boldsymbol{y}%
_{-i}^{T}\boldsymbol{y}_{-i}+\frac{1}{2}y_{i}^{2}\right) \\
=&\overline{G}_{n}\left(\frac{1}{2}\boldsymbol{y}_{-i}^{T}\boldsymbol{y}%
_{-i}+\frac{1}{2}\eta_{\boldsymbol{p},i}^{2}\right)-\overline{G}_{n}\left(%
\frac{1}{2}\boldsymbol{y}_{-i}^{T}\boldsymbol{y}_{-i}+\frac{1}{2}\eta_{%
\boldsymbol{q},i}^{2}\right),
\end{align*}
\begin{align*}
\mathrm{E}[Y_{i}Y_{j}|\boldsymbol{\eta_{p}}\leq\mathbf{Y}\leq\boldsymbol{%
\eta_{q}}]=&\frac{c_{n}}{F_{\mathbf{Y}}(\boldsymbol{\eta_{p}},\boldsymbol{%
\eta_{q}})}\int_{\boldsymbol{\eta}_{\boldsymbol{p},-i}}^{\boldsymbol{\eta}_{%
\boldsymbol{q},-i}}y_{j}\left[\overline{G}_{n}\left(\frac{1}{2}\boldsymbol{y}%
_{-i}^{T}\boldsymbol{y}_{-i}+\frac{1}{2}\eta_{\boldsymbol{p},i}^{2}\right)-%
\overline{G}_{n}\left(\frac{1}{2}\boldsymbol{y}_{-i}^{T}\boldsymbol{y}_{-i}+%
\frac{1}{2}\eta_{\boldsymbol{q},i}^{2}\right)\right]\mathrm{d}\boldsymbol{y}%
_{-i} \\
=&\frac{c_{n}}{F_{\mathbf{Y}}(\boldsymbol{\eta_{p}},\boldsymbol{\eta_{q}})}%
\int_{\boldsymbol{\eta}_{\boldsymbol{p},-ij}}^{\boldsymbol{\eta}_{%
\boldsymbol{q},-ij}}\int_{\eta_{\boldsymbol{p},j}}^{\eta_{\boldsymbol{q}%
,j}}y_{j}\bigg[\overline{G}_{n}\left(\frac{1}{2}\boldsymbol{y}_{-ij}^{T}%
\boldsymbol{y}_{-ij}+\frac{1}{2}y_{j}^{2}+\frac{1}{2}\eta_{\boldsymbol{p}%
,i}^{2}\right) \\
&-\overline{G}_{n}\left(\frac{1}{2}\boldsymbol{y}_{-ij}^{T}\boldsymbol{y}%
_{-ij}+\frac{1}{2}y_{j}^{2}+\frac{1}{2}\eta_{\boldsymbol{q},i}^{2}\right)%
\bigg]\mathrm{d}y_{j}\mathrm{d}\boldsymbol{y}_{-ij},~i\neq j,
\end{align*}
where $\boldsymbol{y}_{-ij}=(y_{1},\cdots,y_{i-1},y_{i+1},%
\cdots,y_{j-1},y_{j+1},\cdots,y_{n})^{T}$.\newline
While
\begin{align*}
&\int_{\eta_{\boldsymbol{p},j}}^{\eta_{\boldsymbol{q},j}}y_{j}\overline{G}%
_{n}\left(\frac{1}{2}\boldsymbol{y}_{-ij}^{T}\boldsymbol{y}_{-ij}+\frac{1}{2}%
y_{j}^{2}+\frac{1}{2}\eta_{\boldsymbol{v},i}^{2}\right)\mathrm{d}%
y_{j}=-\int_{\eta_{\boldsymbol{p},j}}^{\eta_{\boldsymbol{q},j}}\partial_{j}%
\overline{\mathcal{G}}_{n}\left(\frac{1}{2}\boldsymbol{y}_{-ij}^{T}%
\boldsymbol{y}_{-ij}+\frac{1}{2}y_{j}^{2}+\frac{1}{2}\eta_{\boldsymbol{v}%
,i}^{2}\right) \\
&=\overline{\mathcal{G}}_{n}\left(\frac{1}{2}\boldsymbol{y}_{-ij}^{T}%
\boldsymbol{y}_{-ij}+\frac{1}{2}\eta_{\boldsymbol{p},j}^{2}+\frac{1}{2}\eta_{%
\boldsymbol{v},i}^{2}\right)-\overline{\mathcal{G}}_{n}\left(\frac{1}{2}%
\boldsymbol{y}_{-ij}^{T}\boldsymbol{y}_{-ij}+\frac{1}{2}\eta_{\boldsymbol{q}%
,j}^{2}+\frac{1}{2}\eta_{\boldsymbol{v},i}^{2}\right) ,~\boldsymbol{v}\in\{%
\boldsymbol{p,q}\},
\end{align*}
thus
\begin{align*}
\mathrm{E}[Y_{i}Y_{j}|\boldsymbol{\eta_{p}}\leq\mathbf{Y}\leq\boldsymbol{%
\eta_{q}}]=&\frac{1}{F_{\mathbf{Y}}(\boldsymbol{\eta_{p}},\boldsymbol{%
\eta_{q}})}\bigg\{\frac{c_{n}}{c_{n-2,\boldsymbol{p}i,\boldsymbol{p}%
j}^{\ast\ast}}F_{\mathbf{Y}_{\boldsymbol{p}i,\boldsymbol{p}j}}(\boldsymbol{%
\eta}_{\boldsymbol{p},-ij},\boldsymbol{\eta}_{\boldsymbol{q},-ij})-\frac{%
c_{n}}{c_{n-2,\boldsymbol{p}i,\boldsymbol{q}j}^{\ast\ast}}F_{\mathbf{Y}_{%
\boldsymbol{p}i,\boldsymbol{q}j}}(\boldsymbol{\eta}_{\boldsymbol{p},-ij},%
\boldsymbol{\eta}_{\boldsymbol{q},-ij}) \\
&+\frac{c_{n}}{c_{n-2,\boldsymbol{q}i,\boldsymbol{q}j}^{\ast\ast}}F_{\mathbf{%
Y}_{\boldsymbol{q}i,\boldsymbol{q}j}}(\boldsymbol{\eta}_{\boldsymbol{p},-ij},%
\boldsymbol{\eta}_{\boldsymbol{q},-ij}) -\frac{c_{n}}{c_{n-2,\boldsymbol{p}j,%
\boldsymbol{q}i}^{\ast\ast}}F_{\mathbf{Y}_{\boldsymbol{p}j,\boldsymbol{q}i}}(%
\boldsymbol{\eta}_{\boldsymbol{p},-ij},\boldsymbol{\eta}_{\boldsymbol{q}%
,-ij}) \bigg\},~i\neq j.
\end{align*}
In a similar manner, by using integration by parts, after some algebra we
obtain
\begin{align*}
\mathrm{E}[Y_{i}^{2}|\boldsymbol{\eta_{p}}\leq\mathbf{Y}\leq\boldsymbol{%
\eta_{q}}]&=\frac{c_{n}}{F_{\mathbf{Y}}(\boldsymbol{\eta_{p}},\boldsymbol{%
\eta_{q}})}\int_{\boldsymbol{\eta_{p}}}^{\boldsymbol{\eta_{q}}%
}y_{i}^{2}g_{n}\left(\frac{1}{2}\boldsymbol{y}^{T}\boldsymbol{y}\right)%
\mathrm{d}\boldsymbol{y} \\
&=\frac{c_{n}}{F_{\mathbf{Y}}(\boldsymbol{\eta_{p}},\boldsymbol{\eta_{q}})}%
\int_{\boldsymbol{\eta}_{\boldsymbol{p},-i}}^{\boldsymbol{\eta}_{\boldsymbol{%
q},-i}}\int_{\eta_{\boldsymbol{p},i}}^{\eta_{\boldsymbol{q}%
,i}}y_{i}^{2}g_{n}\left(\frac{1}{2}\boldsymbol{y}_{-i}^{T}\boldsymbol{y}%
_{-i}+\frac{1}{2}y_{i}^{2}\right)\mathrm{d}y_{i}\mathrm{d}\boldsymbol{y}_{-i}
\\
&=\frac{c_{n}}{F_{\mathbf{Y}}(\boldsymbol{\eta_{p}},\boldsymbol{\eta_{q}})}%
\int_{\boldsymbol{\eta}_{\boldsymbol{p},-i}}^{\boldsymbol{\eta}_{\boldsymbol{%
q},-i}}\int_{\eta_{\boldsymbol{p},i}}^{\eta_{\boldsymbol{q}%
,i}}-y_{i}\partial_{i}\overline{G}_{n}\left(\frac{1}{2}\boldsymbol{y}%
_{-i}^{T}\boldsymbol{y}_{-i}+\frac{1}{2}y_{i}^{2}\right)\mathrm{d}%
\boldsymbol{y}_{-i} \\
&=\frac{c_{n}}{F_{\mathbf{Y}}(\boldsymbol{\eta_{p}},\boldsymbol{\eta_{q}})}%
\int_{\boldsymbol{\eta}_{\boldsymbol{p},-i}}^{\boldsymbol{\eta}_{\boldsymbol{%
q},-i}} \bigg[\eta_{\boldsymbol{p},i}\overline{G}_{n}\left(\frac{1}{2}%
\boldsymbol{y}_{-i}^{T}\boldsymbol{y}_{-i}+\frac{1}{2}\eta_{\boldsymbol{p}%
,i}^{2}\right)-\eta_{\boldsymbol{q},i}\overline{G}_{n}\left(\frac{1}{2}%
\boldsymbol{y}_{-i}^{T}\boldsymbol{y}_{-i}+\frac{1}{2}\eta_{\boldsymbol{q}%
,i}^{2}\right) \\
&~~~+\int_{\eta_{\boldsymbol{p},i}}^{\eta_{\boldsymbol{q},i}}\overline{G}%
_{n}\left(\frac{1}{2}\boldsymbol{y}_{-i}^{T}\boldsymbol{y}_{-i}+\frac{1}{2}%
y_{i}^{2}\right)\mathrm{d}y_{i}\bigg] \mathrm{d}\boldsymbol{y}_{-i} \\
&=\frac{c_{n}}{F_{\mathbf{Y}}(\boldsymbol{\eta_{p}},\boldsymbol{\eta_{q}})}%
\bigg\{ \frac{\eta_{\boldsymbol{p},i}}{c_{n-1,\boldsymbol{p},i}}F_{\mathbf{Y}%
_{\boldsymbol{p},-i}}(\boldsymbol{\eta}_{\boldsymbol{p},-i},\boldsymbol{\eta}%
_{\boldsymbol{q},-i})-\frac{\eta_{\boldsymbol{q},i}}{c_{n-1,\boldsymbol{q},i}%
}F_{\mathbf{Y}_{\boldsymbol{q},-i}}(\boldsymbol{\eta}_{\boldsymbol{p},-i},%
\boldsymbol{\eta}_{\boldsymbol{q},-i})+\frac{1}{c_{n}^{\ast}}F_{\mathbf{Y}%
^{\ast}}(\boldsymbol{\eta_{p}},\boldsymbol{\eta_{q}})\bigg\}.
\end{align*}
As for $\mathrm{MRVaR}_{(\boldsymbol{\eta_{p}},\boldsymbol{\eta_{q}})}(%
\mathbf{Y})_{k}$, using Theorem 2 we immediately obtain (\ref{(v11)}).
Therefore we obtain $(\ref{(v10)})$, as required.

\noindent $\mathbf{Proof~of~ Proposition~ 3.}$ (i) Using definition of MRVaR
and $\mathbf{Y}=\mathbf{\Sigma}^{-1/2}(\mathrm{ln}\mathbf{Z}-\boldsymbol{\mu}%
)$, we have
\begin{align*}
\mathrm{MRVaR}_{(\boldsymbol{p,q})}(\mathbf{Z})&=\frac{c_{n}}{\sqrt{|\mathbf{%
\Sigma}|}F_{\mathbf{Z}}(\boldsymbol{a,b})}\int_{\boldsymbol{a}}^{\boldsymbol{%
b}}\boldsymbol{z}\left(\prod_{k=1}^{n}z_{k}^{-1}\right)g_{n}\left(\frac{1}{2}%
(\mathrm{ln}\boldsymbol{z-\mu})^{T}\mathbf{\Sigma}^{-1}(\ln\boldsymbol{z-\mu}%
)\right)\mathrm{d}\boldsymbol{z} \\
&=\frac{c_{n}}{F_{\mathbf{Z}}(\boldsymbol{a,b})}\int_{\boldsymbol{\zeta_{p}}%
}^{\boldsymbol{\zeta_{q}}}\mathrm{e}^{\boldsymbol{\mu}+\mathbf{\Sigma}^{1/2}%
\boldsymbol{y}}g_{n}\left(\frac{1}{2}\boldsymbol{y}^{T}\boldsymbol{y}\right)%
\mathrm{d}\boldsymbol{y},
\end{align*}
where
\begin{equation*}
\mathrm{e}^{\boldsymbol{\mu}+\mathbf{\Sigma}^{1/2}\boldsymbol{y}}=\left(%
\mathrm{e}^{\mu_{1}+\mathbf{\Sigma}_{1}^{1/2}\boldsymbol{y}},\mathrm{e}%
^{\mu_{2}+\mathbf{\Sigma}_{2}^{1/2}\boldsymbol{y}},\cdots,\mathrm{e}%
^{\mu_{n}+\mathbf{\Sigma}_{n}^{1/2}\boldsymbol{y}}\right)^{\mathrm{T}},~%
\boldsymbol{a}=\mathrm{VaR}_{\boldsymbol{p}}(\mathbf{Z})~ \mathrm{and} ~%
\boldsymbol{b}=\mathrm{VaR}_{\boldsymbol{q}}(\mathbf{Z}).
\end{equation*}
From Eq. (\ref{(v12)}), for $\forall k\in\{1,2,\cdots,n\}$, we know
\begin{align*}
&\int_{\boldsymbol{\zeta_{p}}}^{\boldsymbol{\zeta_{q}}}\mathrm{e}^{\mu_{k}+%
\mathbf{\Sigma}_{k}^{1/2}\boldsymbol{y}}c_{n}g_{n}\left(\frac{1}{2}%
\boldsymbol{y}^{T}\boldsymbol{y}\right)\mathrm{d}\boldsymbol{y} \\
& =\mathrm{e}^{\mu_{k}}\int_{\boldsymbol{\zeta_{p}}}^{\boldsymbol{\zeta_{q}}%
}\psi_{n}\left(-\frac{1}{2}\sigma_{kk}\right)f_{\mathbf{P}_{k}\dagger}(%
\boldsymbol{y})\mathrm{d}\boldsymbol{y} \\
&=\mathrm{e}^{\mu_{k}}\psi_{n}\left(-\frac{1}{2}\sigma_{kk}\right)F_{\mathbf{%
P}_{k}\dagger}(\boldsymbol{\zeta_{p}},\boldsymbol{\zeta_{q}}),
\end{align*}
thus we obtain (\ref{(v14)}), ending proof of (i).\newline
(ii) By definition of MRCov, we obtain
\begin{align*}
\mathrm{MRCov}_{(\boldsymbol{p},\boldsymbol{q})}(\mathbf{Z})&=\mathrm{E}%
\left[\mathbf{Z}\mathbf{Z}^{T}|\mathrm{VaR}_{\boldsymbol{p}}(\mathbf{Z})\leq%
\mathbf{Z}\leq \mathrm{VaR}_{\boldsymbol{q}}(\mathbf{Z})\right]-\mathrm{MRVaR%
}_{(\boldsymbol{p},\boldsymbol{q})}(\mathbf{Z})\mathrm{MRVaR}_{(\boldsymbol{p%
},\boldsymbol{q})}^{T}(\mathbf{Z}).
\end{align*}
For $\mathrm{E}\left[\mathbf{Z}\mathbf{Z}^{T}|\mathrm{VaR}_{\boldsymbol{p}}(%
\mathbf{Z})\leq\mathbf{Z}\leq\mathrm{VaR}_{\boldsymbol{q}}(\mathbf{Z})\right]
$, let $\mathbf{Y}=\mathbf{\Sigma}^{-1/2}(\mathrm{ln}\mathbf{Z}-\boldsymbol{%
\mu})$, we get
\begin{align*}
\mathrm{E}\left[\mathbf{Z}\mathbf{Z}^{T}|\mathrm{VaR}_{\boldsymbol{p}}(%
\mathbf{Z})\leq\mathbf{Z}\leq\mathrm{VaR}_{\boldsymbol{q}}(\mathbf{Z})\right]%
=\frac{c_{n}}{F_{\mathbf{Z}}(\mathrm{VaR}_{\boldsymbol{p}}(\mathbf{Z}),%
\mathrm{VaR}_{\boldsymbol{q}}(\mathbf{Z}))}\int_{\boldsymbol{\zeta_{p}}}^{%
\boldsymbol{\zeta_{q}}}\mathrm{e}^{\boldsymbol{\mu}+\mathbf{\Sigma}^{1/2}%
\boldsymbol{y}}\left(\mathrm{e}^{\boldsymbol{\mu}+\mathbf{\Sigma}^{1/2}%
\boldsymbol{y}}\right)^{\mathrm{T}}g_{n}\left(\frac{1}{2}\boldsymbol{y}^{T}%
\boldsymbol{y}\right)\mathrm{d}\boldsymbol{y}.
\end{align*}
According to Eq. (\ref{(v13)}), for $\forall k,j\in\{1,2,\cdots,n\}$, we
have
\begin{align*}
&\int_{\boldsymbol{\zeta_{p}}}^{\boldsymbol{\zeta_{q}}}\mathrm{e}^{\mu_{k}+%
\mathbf{\Sigma}_{k}^{1/2}\boldsymbol{y}}\mathrm{e}^{\mu_{j}+\mathbf{\Sigma}%
_{j}^{1/2}\boldsymbol{y}}c_{n}g_{n}\left(\frac{1}{2}\boldsymbol{y}^{T}%
\boldsymbol{y}\right)\mathrm{d}\boldsymbol{y} \\
&=\mathrm{e}^{\mu_{k}+\mu_{j}}\int_{\boldsymbol{\zeta_{p}}}^{\boldsymbol{%
\zeta_{q}}}\mathrm{e}^{(\mathbf{\Sigma}_{k}^{1/2}+\mathbf{\Sigma}_{j}^{1/2})%
\boldsymbol{y}}c_{n}g_{n}\left(\frac{1}{2}\boldsymbol{y}^{T}\boldsymbol{y}%
\right)\mathrm{d}\boldsymbol{y} \\
&=\mathrm{e}^{\mu_{k}+\mu_{j}}\int_{\boldsymbol{\zeta_{p}}}^{\boldsymbol{%
\zeta_{q}}}\psi_{n}\left(-\frac{1}{2}(\sigma_{kk}+2\sigma_{kj}+\sigma_{jj})%
\right)f_{\mathbf{P}_{k,j}\ddagger}(\boldsymbol{y})\mathrm{d}\boldsymbol{y}
\\
&=\mathrm{e}^{\mu_{k}+\mu_{j}}\psi_{n}\left(-\frac{1}{2}(\sigma_{kk}+2%
\sigma_{kj}+\sigma_{jj})\right)F_{\mathbf{P}_{k,j}\ddagger}(\boldsymbol{%
\zeta_{p}},\boldsymbol{\zeta_{q}}).
\end{align*}
Combining (i), and some algebra we obtain (\ref{(v15)}), as required.

\noindent $\mathbf{MRCorrs}$ $\mathbf{of}$ $\mathbf{U}$, $\mathbf{V}$ $%
\mathbf{and}$ $\mathbf{W}$ $\mathbf{for}$ $\mathbf{different}$ $(\boldsymbol{%
p,q})$: Since $\mathrm{MRCorr}$ is a symmetric matrix, we write only its
upper triangular elements for convenience.
\begin{align*}
\mathrm{MRCorr}_{(\boldsymbol{p_{1}},\boldsymbol{q_{1}})}(\mathbf{U})&
=\left(
\begin{array}{ccc}
1 & -0.0008264 & 0.0514972 \\
& 1 & -0.0462871 \\
&  & 1%
\end{array}%
\right) ,~\mathrm{MRCorr}_{(\boldsymbol{p_{2}},\boldsymbol{q_{2}})}(\mathbf{U%
})=\left(
\begin{array}{ccc}
1 & -0.0002787 & 0.0292771 \\
& 1 & -0.0443861 \\
&  & 1%
\end{array}%
\right) , \\
\mathrm{MRCorr}_{(\boldsymbol{p_{3}},\boldsymbol{q_{3}})}(\mathbf{U})&
=\left(
\begin{array}{ccc}
1 & -0.0018052 & 0.0479807 \\
& 1 & -0.0724043 \\
&  & 1%
\end{array}%
\right) ,~\mathrm{MRCorr}_{(\boldsymbol{p_{4}},\boldsymbol{q_{4}})}(\mathbf{U%
})=\left(
\begin{array}{ccc}
1 & -0.0002486 & 0.04172681 \\
& 1 & -0.0360759 \\
&  & 1%
\end{array}%
\right) . \\
\mathrm{MRCorr}_{(\boldsymbol{p_{1}},\boldsymbol{q_{1}})}(\mathbf{V})&
=\left(
\begin{array}{ccc}
1 & 0.3520923 & 0.4203471 \\
& 1 & 0.2921597 \\
&  & 1%
\end{array}%
\right) ,~\mathrm{MRCorr}_{(\boldsymbol{p_{2}},\boldsymbol{q_{2}})}(\mathbf{V%
})=\left(
\begin{array}{ccc}
1 & -0.0030646 & 0.0510202 \\
& 1 & -0.0769316 \\
&  & 1%
\end{array}%
\right) , \\
\mathrm{MRCorr}_{(\boldsymbol{p_{3}},\boldsymbol{q_{3}})}(\mathbf{V})&
=\left(
\begin{array}{ccc}
1 & -0.0046326 & 0.0796897 \\
& 1 & -0.1098094 \\
&  & 1%
\end{array}%
\right) ,~\mathrm{MRCorr}_{(\boldsymbol{p_{4}},\boldsymbol{q_{4}})}(\mathbf{V%
})=\left(
\begin{array}{ccc}
1 & 0.3571721 & 0.4241448 \\
& 1 & 0.2987143 \\
&  & 1%
\end{array}%
\right) . \\
\mathrm{MRCorr}_{(\boldsymbol{p_{1}},\boldsymbol{q_{1}})}(\mathbf{W})&
=\left(
\begin{array}{ccc}
1 & 0.0865982 & 0.1541785 \\
& 1 & 0.02902425 \\
&  & 1%
\end{array}%
\right) ,~\mathrm{MRCorr}_{(\boldsymbol{p_{2}},\boldsymbol{q_{2}})}(\mathbf{W%
})=\left(
\begin{array}{ccc}
1 & -0.0095494 & 0.0710736 \\
& 1 & -0.1056723 \\
&  & 1%
\end{array}%
\right) , \\
\mathrm{MRCorr}_{(\boldsymbol{p_{3}},\boldsymbol{q_{3}})}(\mathbf{W})&
=\left(
\begin{array}{ccc}
1 & -0.0006446 & 0.0992422 \\
& 1 & -0.1151253 \\
&  & 1%
\end{array}%
\right) ,~\mathrm{MRCorr}_{(\boldsymbol{p_{4}},\boldsymbol{q_{4}})}(\mathbf{W%
})=\left(
\begin{array}{ccc}
1 & 0.0743760 & 0.1344204 \\
& 1 & 0.0241352 \\
&  & 1%
\end{array}%
\right) .
\end{align*}

\end{document}